\documentclass[12pt]{amsart}
\usepackage{amsfonts}
\usepackage{latexsym,amsmath,amsthm,amssymb,amsfonts}

\numberwithin{equation}{section}
\allowdisplaybreaks

\def\endproof{$\hfill\Box$\\}
\def\s{\,\,\,\,}

\numberwithin{equation}{section}
\newtheorem{theorem}{Theorem}[section]
\newtheorem{lem}[theorem]{Lemma}
\newtheorem{thm}[theorem]{Theorem}
\newtheorem{pro}[theorem]{Proposition}

\newtheorem{rem}[theorem]{Remark}

\newcounter{Cnumber}

\title[ ]
{\bf Global existence of the solution to Einstein-Maxwell equations with small initial data}
\author[ ]
{Zonglin Jia \qquad Boling Guo }

\date{}
\begin{document}
\maketitle

\begin{abstract}
We study the global existence of Einstein-Maxwell(EM) equations on $\mathbb{R}^4$. We use the method, which relies on wave and Lorentzian gauge conditions, to obtain some exquisite estimates. Our main conclusion is that if the initial data is small enough, then the EM system has a global in time solution.
\end{abstract}

\section{Introduction}
Over the years, the relativistic theories of gravitation and electromagnetism have experienced an intensive and rapidly developing research activity. The system that unifies these two theories consists of the EM equations. Most researches on EM equations are performed under special metrics. However, we study this system in a general metric and get global existence with small initial data.

Before stating our main result, let us explain the problem to be solved. First of all, we make the following agreement. Throughout the paper, the same indices appearing twice means summing it. Besides, we also appoint that, when denoting superscripts or  subscripts without special statements, the Greek letters such as $\alpha,\beta,\gamma,\cdots$ belong to $\{0,1,2,\cdots,n\}$, while the Latin letters $i,j,k,\cdots$ are in $\{1,2,\cdots,n\}$.

A Lorentzian manifold $(V,g)$ is an $n+1$-dimensional smooth manifold $V$ admitting a symmetric 2-tensor $g$ with signature $(-,\underbrace{+,\cdots,+}\limits_{\mbox{$n$-times}})$. In this paper, we consider the following Einstein-Maxwell(EM) equations on $(V,g)$
\begin{eqnarray*}
\left\{
\begin{array}{rcl}
&Ric_{\alpha\beta}-\frac{1}{2}R\cdot g_{\alpha\beta}=T_{\alpha\beta}        & \quad\quad\text{(Einstein equation)}\\
&d\tilde{F}=0      & \quad\quad\text{(the first set of Maxwell equation)}\\
&div_{g}\tilde{F}=0,     & \quad\quad\text{(the second set of Maxwell equation)}
\end{array} \right.
\end{eqnarray*}
where $Ric$ and $R$ are Ricci tensor and scalar curvature of $g$ respectively. Let $(x^{\alpha})$ be a local coordinate chart of $V$. Sometimes, we write $\frac{\partial}{\partial x^{\alpha}}$ as $\partial_{\alpha}$. Furthermore, $\tilde{F}:=\frac{1}{2}\tilde{F}_{\alpha\beta}dx^{\alpha}\wedge dx^{\beta}$(we have assume that $\tilde{F}_{\alpha\beta}=-\tilde{F}_{\beta\alpha}$) is a 2-form and $T_{\alpha\beta}$ is given by
\begin{eqnarray}\label{9}
T_{\alpha\beta}:=\tilde{F}_{\alpha}^{\lambda}\tilde{F}_{\beta\lambda}-\frac{1}{4}g_{\alpha\beta}\tilde{F}^{\lambda\mu}\tilde{F}_{\lambda\mu},
\end{eqnarray}
where we have raised the indices by $(g^{\alpha\beta})$ which is the inverse of the metric matrix $(g_{\alpha\beta})$. From (10.2) of Chapter 6 of \cite{CB} it follows that the Einstein equation is equivalent to a simpler one
\begin{eqnarray}\label{0}
Ric_{\alpha\beta}=\tilde{F}_{\alpha}^{\lambda}\tilde{F}_{\beta\lambda}-\frac{1}{2(n-1)}g_{\alpha\beta}\tilde{F}^{\lambda\mu}\tilde{F}_{\lambda\mu}.
\end{eqnarray}

This paper mainly focuses on EM equations on $\mathbb{R}^3\times\mathbb{R}\equiv\mathbb{R}^4$, which is a simply connected manifold. Hence the first set of Maxwell equation is equivalent to $\tilde{F}=dA$ for some one-form $A$.

Let $h:=g-m$, where $m$ is the Minkowski metric of $\mathbb{R}^4$, and set
\begin{eqnarray*}
\mathfrak{L}:=\{\partial_{\alpha},\s\Omega_{\alpha\beta}:=-x_{\alpha}\partial_{\beta}+x_{\beta}\partial_{\alpha},\s \tilde{S}:=t\partial_t+r\partial_r\s|\alpha,\beta=0,1,2,3\}.
\end{eqnarray*}
This family of vector fields plays a critical role in the study of the wave equation in Minkowski space-time. We denote the above vector fields by $Z^{\iota}$ with an 11-dimensional integer index $\iota=(0,\cdots,1,\cdots,0)$. Let $I:=(\iota_1,\cdots,\iota_k)$, where $|\iota_i|=1$ for $1\leqslant i\leqslant k$, be a multi-index of length $|I|=k$ and let $Z^I:=Z^{\iota_1}\cdots Z^{\iota_k}$ denote a product of $k$ vector fields from the family $\mathfrak{L}$. By a sum $I_1+I_2=I$ we mean a sum over all possible order preserving partitions of the multi-index $I$ into two multi-indices $I_1$ and $I_2$, i.e. if $I=(\iota_1,\cdots,\iota_k)$, then $I_1=(\iota_{i_1},\cdots,\iota_{i_n})$ and $I_2=(\iota_{i_{n+1}},\cdots,\iota_{i_k})$, where $i_1,\cdots,i_k$ is any reordering of the integers $1,\cdots,k$ such that $i_1<\cdots<i_n$ and $i_{n+1}<\cdots<i_k$. With this convention the Leibnitz rule becomes
\begin{eqnarray*}
Z^I(f\cdot g)=\sum\limits_{I_1+I_2=I}(Z^{I_1}f)\cdot(Z^{I_2}g).
\end{eqnarray*}
Let $\Box:=m^{\alpha\beta}\partial_{\alpha}\partial_{\beta}$(where $\alpha,\beta\in\{0,1,2,3\}$). From Section 2 of \cite{LR2010} it follows that we have the next commutation properties:
\begin{eqnarray*}
[\Box,\partial_{\alpha}]=0,\s\s\s\s[\Box,\Omega_{\alpha\beta}]=0\s\s\s\s\mbox{and}\s\s\s\s[\Box,\tilde{S}]=2\Box,
\end{eqnarray*}
where $[X,Y]:=XY-YX$ is the commutator. For $Z\in\mathfrak{L}$, denote $[Z,\Box]:=c_Z\cdot\Box$, i.e. $c_Z=2$, if $Z=\tilde{S}$, and 0, otherwise. Moreover, set
\begin{eqnarray}\label{95}
\mathcal{E}_N(t):=\sum\limits_{|I|\leqslant N,Z\in\mathfrak{L}}\left(||\sqrt{w}\partial Z^Ih(\cdot,t)||_{L^2}+||\sqrt{w}\partial Z^IA(\cdot,t)||_{L^2}\right),
\end{eqnarray}
where
\begin{eqnarray*}
w(q):=
\left\{
\begin{array}{rcl}
&1+(1+|q|)^{1+2\gamma}        & \quad\quad\text{when $q>0$,}\\
&1+(1+|q|)^{-2\mu}         & \quad\quad\text{when $q<0$.}
\end{array} \right.
\end{eqnarray*}
with $q:=|x|-t$ and two constants $\gamma\in(0,1/2),\mu\in(0,1/2)$ being fixed.\\

Now we state the main result of this paper.
\begin{thm}\label{thm6.1}
Given an integer $N\geqslant4$, there exists a constant $\varepsilon_0>0$ such that if $\varepsilon\leqslant\varepsilon_0$ and the initial datum $h|_{t=0}$, $\partial_th|_{t=0}$, $A|_{t=0}$ and $\partial_tA|_{t=0}$ obey $\mathcal{E}_N(0)\leqslant\varepsilon$, then the solution of EM equations $(g(t)=h(t)+m,A(t))$ belongs to the space $\tilde{E}_{N+1}(\infty)$.
\end{thm}

\begin{rem}
The definition of $\tilde{E}_{N+1}(\infty)$ is given in Subsection \ref{subsection2.7}. Moreover, at this time we take $M:=\mathbb{R}^3$ and let $e$ be the standard Euclidean metric in that subsection.
\end{rem}

Let us recall some previous results about Einstein equations. It is well known that these equations are quasilinear wave equations in wave coordinates. In general, to ensure such a system to have global solutions for small initial data there must be ``null condition''. Unfortunately, the condition does not hold for Einstein equations. Hence, in \cite{LR2003}, H. Lindblad and I. Rodnianski brought in the notion of ``weak null condition" and demonstrated that the Einstein equations in harmonic gauges satisfy this condition. In \cite{LR2005}, they used the weak null condition to prove the global existence for the Einstein vacuum equations in wave coordinates, which contradicts beliefs that wave coordinates are ``unstable in the large". The outstanding global problem, which for a long time remained open, was firstly solved by D. Christodoulou and S. Klainerman in \cite{CK1993}. In \cite{LR2010}, Lindblad and Rodnianski used their ingenious approach to show that Einstein-scalar field equations admit global stability. And the idea of our present article does come from \cite{LR2010}. Applying different methods, in \cite{Zip2000} N. Zipser got the conclusion that the trivial solution of Einstein-Maxwell equations admits global nonlinear stability. His way is not to use wave gauge conditions. Instead, he considered generalized energy inequalities associated with Bel-Robinson energy-momentum tensor, designed to mimic the rotation and the conformal Morawetz vector fields of the Minkowski space-time.

Except for the above cases, people also care about Einstein equations coupling with other systems of physical significance in special metrics. In \cite{ST2014}, M. Sango and C. Tadmon analyze the Einstein-Maxwell-Euler equations for an irrotational stiff fluid. Under the spherical symmetry assumption on the space-time, in Bondi coordinates, the considered model is reduced to a nonlinear evolution system of partial integrodifferential equations. In \cite{C2003}, D. Chae proved the global unique existence of classical solutions to the Einstein-Maxwell-Higgs(EMH) systems for small initial data under spherical symmetry. He also obtained the decay estimates of the solutions and found that the corresponding space-time is time-like and null geodesically complete toward the future. In \cite{TT2012}, C. Tadmon and S.B. Tchapnda revisited and generalized, to the Einstein-Yang-Mills-Higgs system, previous results of Christodoulou and Chae concerning global solutions for the Einstein-scalar field and the EMH equations. In \cite{D2003}, M. Dafermos considered a trapped characteristic initial value problem for the spherically symmetric Einstein-Maxwell-scalar field equations. For an open set of initial data whose closure contains in particular Reissner-Nordstr\"om data, the future boundary of the maximal domain of development is found to be a light-like surface along which the curvature blow up, and yet the metric can be continuously extended beyond it. His result is related to the strong cosmic censorship conjecture of Roger Penrose. In \cite{SW1997}, J.A. Smoller and A.G. Wasserman proved that any solution to the spherically symmetric $SU(2)$ Einstein-Yang-Mills equations, which is defined in the far field and is asymptotically flat, is globally defined. Their result applies in particular to the interior of colored black holes.

Furthermore, the Maxwell equation in non-relativistic situations is very important. In \cite{Y2015}, Yang presented a new approach to study the asymptotic behavior of the global solutions of the Maxwell-Klein-Gordon(MKG) equations. This is the first result to give a complete and precise description of such behavior of large nonlinear charged scalar fields. And it is conjectured that the solutions of the MKG equations disperse as linear waves and enjoy peeling properties for pointwise estimates. In \cite{YY2018}, Yang and Yu provided a gauge independent proof of this conjecture.

Now let us briefly introduce the method we use. At the beginning, in the process of getting local well-posedness, we brought in wave and Lorentzian gauges to transform the EM equations into hyperbolic systems called the reduced EM systems. In order to show that the solution to the reduced EM systems also solves the full set of EM equations, we have to require that the initial datum sets satisfy Einstein-Maxwell constraints. For the details of the above concepts readers may refer to Section \ref{section2} of this paper. Secondly, as soon as we get a local solution, it is natural to consider the maximal existence time $T$ and assume it to be finite. In the next, we define $T^*$ to be
\begin{eqnarray}\label{00}
T^*:=\sup\{T_0:\exists C=C(T),\,s.t.\s\forall t\in[0,T_0),\s\mathcal{E}_N(t)\leqslant2C\varepsilon\}
\end{eqnarray}
and suppose that $T^*<T$. We will show that if $\varepsilon>0$ is small enough, then the inequality in (\ref{00}) implies the same inequality with $2C$ replaced by $C$ for all $t<T^*$. This contradicts the maximality of $T^*$ and it follows that the inequality holds for all $t<T$. Moreover, since the energy $\mathcal{E}_N(t)$ is now finite at $t=T$, we can extend the solution beyond $T$ thus contradicting maximality of $T$ and showing that $T=\infty$. Thirdly, we must get several decay estimates to ensure the smooth implementation of the above process. The critical part is Theorem \ref{thm6.2} and this involves some techniques of controlling the inhomogeneous terms on the right hand side of the reduced EM equations. At last, our task is to get energy estimates. Noting the definition of $\mathcal{E}_N$, we compute $\overset{\sim}{\Box}_gZ^Ih$ and $\overset{\sim}{\Box}_gZ^IA$ where $\overset{\sim}{\Box}_g:=g^{\alpha\beta}\partial_{\alpha}\partial_{\beta}$ and $\alpha,\beta\in\{0,1,2,3\}$. Applying Proposition 6.2 of \cite{LR2010} and Gronwall inequality leads to the needed results.

This paper is organized as follows. We devote Section \ref{section2} to some notations and preliminaries. In Section \ref{section3}, we reduce the EM equations to a hyperbolic system and obtain the local well-posedness. In Section \ref{section4}, we focus on the EM system of $(h=g-m,A)$ on $\mathbb{R}^4$. Finally, the decay estimates and energy estimates are presented in Section \ref{section5} and Section \ref{section6} respectively.

\section{Notations and Preliminaries}\label{section2}
In this article, the symbol ``$Q_1\lesssim Q_2$" means that there exists a constant $C$ such that $Q_1\leqslant C\cdot Q_2$ for two given quantities $Q_1$ and $Q_2$.  Throughout our paper, the constant $C$ may depend upon the maximal existence time $T$.
\subsection{Wave gauges}
Let $W$ be a smooth manifold and $\hat{e}$ and $g$ are two Riemannian or pseudo Riemannian metric on $W$. We say that $g$ is in wave gauge with respect to $\hat{e}$ if the identity map $id:W\longrightarrow W$ is a wave map from $(W,g)$ into $(W,\hat{e})$, i.e.
\begin{eqnarray*}
\hat{F}:=\mbox{tr}_g(\Gamma-\hat{\Gamma})\equiv g^{\alpha\beta}(\Gamma^{\lambda}_{\alpha\beta}-\hat{\Gamma}^{\lambda}_{\alpha\beta})\partial_{\lambda}=0,
\end{eqnarray*}
where $\Gamma^{\lambda}_{\alpha\beta}$ and $\hat{\Gamma}^{\lambda}_{\alpha\beta}$ are the connection coefficients of $g$ and $\hat{e}$ respectively.
\subsection{Sobolev space on Riemannian manifold}
Given a smooth Riemannian manifold $(M,e)$, the Sobolev space $W_s^p(M,e)$ or $W_s^p$ for short($1\leqslant p\leqslant\infty$) is the space of functions, or tensor fields of some given type, on $(M,e)$ with $L^p$ integrable generalized derivatives of order not bigger than $s$ in
the metric $e$. More precisely, it is a Banach space with the norm
\begin{eqnarray}\label{norm1}
||f||^p_{W^p_s}:=\sum\limits_{0\leqslant k\leqslant s}\int_M|\tilde{D}^k f|^p\,\mu_e,\s\s\s\s\mbox{for}\s\s 1\leqslant p<\infty;
\end{eqnarray}
and
\begin{eqnarray}\label{norm2}
||f||_{W^{\infty}_s}:=\max\limits_{0\leqslant k\leqslant s}||\tilde{D}^k f||_{L^{\infty}}
\end{eqnarray}
with $|\cdot|$ the pointwise norm of tensors in the metric $e$, where $\tilde{D}$ is the covariant derivative and $\mu_e$ is volume element with respect to $e$.

The Sobolev spaces $\overset{0}{W_s^p}\subseteq W_s^p$ on $(M,e)$ are defined as the closures with respect to the norms (\ref{norm1}) or (\ref{norm2}) of spaces of $C^{\infty}$-functions, or tensor fields of some given type, with compact support in $M$. For simplicity, we usually denote $W^2_s$ by $H_s$ and $\overset{0}{W^2_s}$ by $\overset{0}{H_s}$.

We use $\bar{C}^k$ to denote the space of functions, or tensor fields of some given type, which admit $k$-derivatives and the $k$-th derivative is continuous and bounded.
\subsection{Sobolev regularity}
We say that the smooth Riemannian manifold $(M,e)$ is Sobolev regular if the spaces $W_s^p$ defined on $(M,e)$ satisfy the embedding and multiplication properties. For the details of the above two properties the readers can refer to Proposition 2.1 and 2.3 of Appendix 1 of \cite{CB}.
\begin{thm}\label{thm2.1}
Complete Riemannian manifolds, or Riemannian manifolds with Lipshitzian boundary, are Sobolev regular.
\end{thm}
\textbf{Proof.} We refer to Theorem 2.5 of Appendix 1 of \cite{CB} for the details.
\endproof
\subsection{Uniform equivalence of Riemannian metric}
Given two Riemannian metric $e_1$ and $e_2$ on a smooth manifold $M$, we say they are uniformly equivalent to each other if there exist two constants $0<A<B<\infty$ such that for all $X\in TM$, the following holds true
\begin{eqnarray*}
A\cdot e_1(X,X)\leqslant e_2(X,X)\leqslant B\cdot e_1(X,X).
\end{eqnarray*}

\subsection{Sliced spacetime}
A Lorentzian spacetime $(V_T,g)$ with $V_T:=M\times(-T,T)$ is called sliced if there exists a time-dependent vector $\beta(x,t):=\beta^i(x,t)\frac{\partial}{\partial x^i}(x)$, which is called the shift, tangent to the space slice $M_t:=M\times\{t\}$($(x^i)$ is a local coordinate chart of $M$) such that

(1) $\vec{n}(\cdot,t):=\frac{\partial}{\partial t}(\cdot,t)-\beta(\cdot,t)$ is a normal vector to $M_t$ for all $t\in(-T,T)$. That is to say, given $X_t\in T(M_t)$($T(M_t)$ is the tangent bundle of $M_t$), we have
\begin{eqnarray}\label{5}
g(\vec{n}(\cdot,t),X_t)=0.
\end{eqnarray}
It is easy to check that (\ref{5}) is equivalent to $g_{0i}=g_{ij}\beta^j$;

(2) $\vec{n}$ is timelike, namely,
\begin{eqnarray}\label{6}
g(\vec{n},\vec{n})<0.
\end{eqnarray}
Combining (\ref{5}) and (\ref{6}) we arrive at $g_{ij}\beta^i\beta^j>g_{00}$.

Using the above inequality we define a positive function $N$ called the lapse which is given by $N:=\sqrt{g_{ij}\beta^i\beta^j-g_{00}}$. Then we can write $g$ as
\begin{eqnarray}\label{11}
g=-N^2dt\otimes dt+g_{ij}(dx^i+\beta^idt)\otimes(dx^j+\beta^jdt).
\end{eqnarray}
Because $g$ is Lorentzian, $(g_{ij})$ is positive definite. In order to represent $(g^{\alpha\beta})$ via $(N,\beta,g_{ij})$, we denote the inverse of $(g_{ij})$ by $(g_{*}^{ij})$. It is easy to check
\begin{eqnarray*}
g^{00}=-N^{-2},\s\s\s\s g^{0i}=N^{-2}\beta^i\s\s\s\s\mbox{and}\s\s\s\s g^{ij}=g_*^{ij}-N^{-2}\beta^i\beta^j.
\end{eqnarray*}

Thanks to the above discussion, now we can give the following lemma without proof.
\begin{lem}\label{lem2.3}
A spacetime $(V_T,g)$ is sliced and Lorentzian if and only if $g_t:=i_t^*g$, which is induced by the embedding $i_t:M\longrightarrow V_T$, $x\mapsto(x,t)$, is positive definite and
\begin{eqnarray*}
g_{0i}\cdot g_{*}^{ij}\cdot g_{0j}>g_{00}.
\end{eqnarray*}
\end{lem}

\begin{rem}
Lemma \ref{lem2.3} tells us that determining a sliced Lorentzian metric $g$ on $V_T$ is equivalent to determining the following quantities:

(1) a Riemannian metric $g_t$ on $M_t$;

(2) a positive function $N$ on $M_t$;

(3) a tangent vector field $\beta$ to $M_t$.
\end{rem}

\begin{rem}
In case $g_t$ is positive definite, the inverse of $(g_t)_{ij}$ which is denoted by $(g_t)^{ij}$ is just $g_*^{ij}$. Hence, throughout this article we always use the symbol $g_*^{ij}$.
\end{rem}

\subsection{Regular sliced spacetime}
Given a Riemannian metric $e$ on $M$, a sliced spacetime $(V_T,g,\beta)$ is called regular with respect to $e$ if

(1) The metrics $g_t$ are uniformly equivalent to $e$, i.e. there exist continuous strictly positive functions $B_1(t)$, $B_2(t)$ such that for all $t\in(-T, T)$ and each tangent vector $X$ to $M$ it holds true on $M$
\begin{eqnarray*}
B_1(t)\cdot e(X,X)\leqslant g_t(X,X)\leqslant B_2(t)\cdot e(X,X);
\end{eqnarray*}

(2) The lapse $N$ is such that there exist continuous strictly positive functions $C_1(t)$, $C_2(t)$ on $(-T, T)$ such that on each $M_t$ it holds true
\begin{eqnarray*}
C_1(t)\leqslant N(x,t)\leqslant C_2(t)\s\s\s\s\forall x\in M;
\end{eqnarray*}

(3) The shift $\beta$ is uniformly bounded in $e$-norm on each $M_t$ by a number $b(t)$.

\subsection{Sobolev space on $V_T$}\label{subsection2.7}
We denote by $C_0^s(T)$ the restriction to $V_T$ of $C^s$-functions or tensor fields with compact support in $V:=M\times\mathbb{R}$ and by $E_s(T)$ the following Banach space
$$E_s(T):=\bigcap\limits_{0\leqslant k\leqslant s}C^{s-k}((-T,T),H_k(M,e)).$$
Besides, the next Banach space is important: $$\overset{0}{E_s}(T):=\bigcap\limits_{0\leqslant k\leqslant s}C^{s-k}((-T,T),\overset{0}{H_k}(M,e)).$$
It is the completion of $C_0^s(T)$ in the $E_s(T)$ norm:
\begin{eqnarray*}
||u||_{E_s(T)}:=\sup\limits_{t\in(-T,T)}\max\limits_{0\leqslant k\leqslant s}||(\partial_t^{s-k}u)(\cdot,t)||_{H_k(M,e)}.
\end{eqnarray*}
For more details we refer to Definition 2.23 of Appendix 3 of \cite{CB}.

The Sobolev space $\tilde{E}_s(T)$ is the space of functions or tensor fields $u$, such that $u\in\bar{C}(V_T)$, space of continuous and bounded tensor fields or functions on $V_T$, while $\partial_tu, \tilde{D}u\in\overset{0}{E}_{s-1}(T)$, where $\tilde{D}$ is the Levi-Civita connection of $(M,e)$.
\subsection{Einstein-Maxwell initial data set}\label{sub2.8}
An Einstein-Maxwell initial data set is a  seven-tuple $(M;\bar{g},K,\bar{F},\bar{E},\bar{\beta},\bar{N})$, where $(M,\bar{g})$ is an $n$-dimensional Riemannian manifold with Riemannian metric $\bar{g}$. $K$ is a symmetric 2-tensor on $M$. Meanwhile, $\bar{F}$ is a 2-form and $\bar{E}$ and $\bar{\beta}$ are two vector fields on $M$. Moreover, $\bar{N}$ is a positive function on $M$.
\subsection{Einstein-Maxwell development}
The development of the initial data set $(M;\bar{g},K,\bar{F})$ is a triplet $(V_T,g,\tilde{F})$ with $V_T:=M\times(-T,T)$ for some $T>0$, a Lorentzian metric $g$ and a 2-form $\tilde{F}$ on $V_T$, such that the embedding $i_0$ of $M$ into $V_T$(Recall that $i_0(x):=(x,0)$ for all $x\in M$) enjoys the following properties:

(a) The metric $\bar{g}$ is the pullback of $g$ by $i_0$, i.e. $\bar{g}=i_0^*g$;

(b) $K$ is the second fundamental form of $i_0(M)$ as a submanifold $(V_T,g)$;

(c) The 2-form $\bar{F}$ is the pullback of $\tilde{F}$ by $i_0$, i.e. $\bar{F}=i_0^*\tilde{F}$.

Furthermore, $(g,\tilde{F})$ satisfies on $V_T$ the Einstein-Maxwell equations and $(V_T,g)$ is a sliced spacetime.
\subsection{Einstein-Maxwell constraints}
Restricting EM equations to the initial data set $(M;\bar{g},K,\bar{F})$, which is embedded into $V_T$, leads to the following identities called the constraints
\begin{eqnarray}\label{1}
\bar{R}-|K|^2_{\bar{g}}+(\mbox{tr}_{\bar{g}}K)^2=2\bar{N}^{-2}\cdot T(\partial_t-\beta,\partial_t-\beta)|_{t=0},
\end{eqnarray}

\begin{eqnarray}\label{2}
(\mbox{div}_{\bar{g}}K)_i-\partial_i(\mbox{tr}_{\bar{g}}K)=-\bar{N}^{-1}\cdot T(\partial_t-\beta,\partial_i)|_{t=0},
\end{eqnarray}

\begin{eqnarray}\label{3}
\bar{\nabla}_h\bar{F}_{ij}+\bar{\nabla}_j\bar{F}_{hi}+\bar{\nabla}_i\bar{F}_{jh}=0,
\end{eqnarray}

\begin{eqnarray}\label{10}
\tilde{F}_{0i}|_{t=0}:=\bar{N}\cdot\bar{E}^j\cdot\bar{g}_{ji}
\end{eqnarray}
and
\begin{eqnarray}\label{4}
\bar{\nabla}_i\bar{E}^i=0,
\end{eqnarray}
where $\bar{\nabla}$ and $\bar{R}$ are the Levi-Civita connection and the scalar curvature of $(M,\bar{g})$ respectively. The details of the above constraints can be found in Section 4.1 and 10.1 of Chapter 6 in \cite{CB}.
\begin{rem}
Unless we give (\ref{10}), one can not determine the values on the right hand side of (\ref{1}) and (\ref{2}). In the next, we are going to write their specific expressions via $(\bar{g},\bar{N},\bar{F},\bar{E},\bar{\beta})$:
\begin{eqnarray*}
2\bar{N}^{-2}\cdot T(\partial_t-\beta,\partial_t-\beta)|_{t=0}&=&\bar{g}^{ij}(\bar{\beta}^k\bar{F}_{ki}-\bar{N}\bar{E}^l\bar{g}_{li})(\bar{\beta}^k\bar{F}_{kj}-\bar{N}\bar{E}^l\bar{g}_{lj})/\bar{N}^2\\
&&+\bar{g}^{ik}\bar{g}^{jl}\bar{F}_{ij}\bar{F}_{kl}/2+2\bar{g}_{ij}\bar{E}^i\bar{E}^j
\end{eqnarray*}
and
\begin{eqnarray*}
-\bar{N}^{-1}\cdot T(\partial_t-\beta,\partial_i)|_{t=0}=\bar{E}^l\bar{F}_{li}+\bar{N}^{-1}\bar{\beta}^k\bar{g}^{pq}\bar{F}_{pk}\bar{F}_{qi}.
\end{eqnarray*}
\end{rem}
\subsection{The null frame}\label{sub2.11}
At each point $x\in\mathbb{R}^4$, we introduce a pair of null vectors $(L,\underline{L})$ with
$$L:=\partial_t+\partial_r\s\s\s\s\mbox{and}\s\s\s\s\underline{L}:=\partial_t-\partial_r,$$
where $x:=(x^0,x^1,x^2,x^3)$, $t:=x^0$ and $r:=\sqrt{(x^1)^2+(x^2)^2+(x^3)^2}$. Sometimes, we also denote $x^{i}$ by $x_{i}$($i=1,2,3$) and $L$ by $\bar{\partial}_0$. Let $S_1$ and $S_2$ be two orthonormal smooth tangent vector fields to the unit sphere $\mathbb{S}^2$, where the orthogonality is in the sense of the standard  metric of $\mathbb{S}^2$. For convenience $A,B,C,D,\cdots$ means any of the vectors $S_1$ and $S_2$ at times. Given a 1-tensor $\pi:=\pi_{\beta}dx^{\beta}$ and a 2-tensor $p:=p_{\alpha\beta}dx^{\alpha}\otimes dx^{\beta}$($\alpha,\beta\in\{0,1,2,3\}$), we define $\pi_X:=\pi(X)$ and $p_{XY}:=p(X,Y)$, provided $X, Y$ are two vector fields.

Suppose that
$$eu:=dt\otimes dt+\sum\limits_{i=1}^3dx^i\otimes dx^i$$
is the standard Euclidean metric of $\mathbb{R}^4$. Then we have
\begin{eqnarray*}
eu_{\underline{L}L}=eu_{LA}=eu_{\underline{L}A}=0,\s\s\s\s eu_{LL}=eu_{\underline{L}\underline{L}}=2
\end{eqnarray*}
and
\begin{eqnarray*}
eu_{AB}=\delta_{AB}:=\left\{
\begin{array}{rcl}
&0      & A\not=B\\
&1      & A=B.
\end{array} \right.
\end{eqnarray*}
The inverse of $eu$ is
\begin{eqnarray*}
eu^{\underline{L}L}=eu^{LA}=eu^{\underline{L}A}=0,\s\s\s\s eu^{LL}=eu^{\underline{L}\underline{L}}=1/2,\s\s\s\s eu^{AB}=\delta_{AB}.
\end{eqnarray*}

Noting that $\tilde{S}_1$ and $\tilde{S}_2$ are defined only locally on $\mathbb{S}^2$, we replace them with the projections
\begin{eqnarray*}
\bar{\partial}_i:=\partial_i-\omega_i\cdot\partial_r,\s\s\s\s \omega_i:=x_i/r\s\s\mbox{and}\s\s i=1,2,3.
\end{eqnarray*}
It is nor hard to see that $\{\bar{\partial}_1,\bar{\partial}_2,\bar{\partial}_3\}$ gives a set of global and linear dependent vector fields on $\mathbb{S}^2$. Moreover, one can also represent $\{\bar{\partial}_i|i=1,2,3\}$ by $\tilde{S}_1$ and $\tilde{S}_2$, i.e.
\begin{eqnarray}\label{19}
\bar{\partial}_i=\tilde{S}_1^i\cdot \tilde{S}_1+\tilde{S}_2^i\cdot \tilde{S}_2,
\end{eqnarray}
where $\tilde{S}_j:=\tilde{S}_j^i\cdot\partial_i$ and $j=1,2$.

We call $\{L,\underline{L},\tilde{S}_1,\tilde{S}_2\}$ the null frame and introduce the coming notation. Let $\mathcal{T}:=\{L,\tilde{S}_1,\tilde{S}_2\}$, $\mathcal{U}:=\{\underline{L},L,\tilde{S}_1,\tilde{S}_2\}$, $\mathcal{L}:=\{L\}$ and $\mathcal{S}:=\{\tilde{S}_1,\tilde{S}_2\}$. For any $l$ of these families $\mathcal{V}_1,\cdots,\mathcal{V}_l$(namely, $\mathcal{V}_1,\cdots,\mathcal{V}_l\in\{\mathcal{T},\mathcal{U},\mathcal{L},\mathcal{S}\}$) and an arbitrary $k$-tensor $p:=p_{\alpha_1\cdots\alpha_k}dx^{\alpha_1}\otimes\cdots\otimes dx^{\alpha_k}$($\alpha_j\in\{0,1,2,3\}$ and $1\leqslant j\leqslant k$) with $k\geqslant l$, we define two norms $|p|$ and $|p|_{\mathcal{V}_1\cdots\mathcal{V}_l}$ as
\begin{eqnarray*}
|p|^2:=\sum\limits_{\alpha_1,\cdots,\alpha_k=0}^3(p_{\alpha_1\cdots\alpha_k})^2
\end{eqnarray*}
and
\begin{eqnarray*}
|p|^2_{\mathcal{V}_1\cdots\mathcal{V}_l}:&=&\sum\limits_{V_1,V'_1\in\mathcal{V}_1}\cdots\sum\limits_{V_l,V'_l\in\mathcal{V}_l}\sum\limits_{\alpha_{l+1},\cdots,\alpha_k=0}^3eu^{V_1V_1'}\cdots eu^{V_lV'_l}\\
&&\cdot p(V_1,\cdots,V_l,\partial_{\alpha_{l+1}},\cdots,\partial_{\alpha_k})\cdot p(V'_1,\cdots,V'_l,\partial_{\alpha_{l+1}},\cdots,\partial_{\alpha_k})
\end{eqnarray*}
It is not difficult to check that $|p|_{\mathcal{V}_1\cdots\mathcal{V}_l}$ is independent of the choice of $\{\tilde{S}_1,\tilde{S}_2\}$ on $\mathbb{S}^2$.
\subsection{The Minkowski metric}
The Minkowski metric $m$ of $\mathbb{R}^4$ is given by
\begin{eqnarray*}
m:=-dt\otimes dt+\sum\limits_{i=1}^3dx^i\otimes dx^i.
\end{eqnarray*}

From Section 4 of \cite{LR2010} it follows that
\begin{eqnarray*}
m_{LL}=m_{\underline{L}\underline{L}}=m_{LA}=m_{\underline{L}A}=0,\s\s\s\s m_{L\underline{L}}=m_{\underline{L}L}=-2\s\s\s\s m_{AB}=\delta_{AB}
\end{eqnarray*}
The inverse of the metric has the form
\begin{eqnarray*}
m^{LL}=m^{\underline{L}\underline{L}}=m^{LA}=m^{\underline{L}A}=0,\s\s\s\s m^{L\underline{L}}=m^{\underline{L}L}=-1/2,\s\s\s\s m^{AB}=\delta^{AB}.
\end{eqnarray*}

Let $\partial$ be the Levi-Civita connection of $(\mathbb{R}^4,m)$. We shall use it to define a new differential operator $\bar{\partial}$ as follows. Provided $p$ is a $k$-tensor and $q:=r-t$, $\bar{\partial}p$ is given by
\begin{eqnarray*}
\bar{\partial}p:=\partial p-\partial_{\partial r}p\otimes dq,
\end{eqnarray*}
where we recall $\partial_r:=\frac{\partial}{\partial r}$. Easily, the readers, reviewing the definition of $\{\bar{\partial}_{\beta}|\beta=0,1,2,3\}$ in Subsection \ref{sub2.11}, can check that
\begin{eqnarray*}
\bar{\partial}p=\bar{\partial}_{\beta}\left(p_{\alpha_1\cdots\alpha_k}\right)\cdot dx^{\alpha_1}\otimes\cdots\otimes dx^{\alpha_k}\otimes dx^{\beta}.
\end{eqnarray*}

\begin{lem}\label{lem2.7}
At a point $(x,t)\in\mathbb{R}^4$, we have
\begin{eqnarray*}
|\bar{\partial}p|^2=|\partial_Lp|^2+|\partial_{\tilde{S}_1}p|^2+|\partial_{\tilde{S}_2}p|^2.
\end{eqnarray*}
\end{lem}
\textbf{Proof.} By the definitions, the following calculation is trival
\begin{eqnarray*}
|\bar{\partial}p|^2&=&\sum\limits_{\alpha_1,\cdots,\alpha_k,\beta}|\bar{\partial}_{\beta}(p_{\alpha_1\cdots\alpha_k})|^2\\
&=&\sum\limits_{\alpha_1,\cdots,\alpha_k}|L(p_{\alpha_1\cdots\alpha_k})|^2+\sum\limits_{\alpha_1,\cdots,\alpha_k}\sum\limits_{i=1}^3|\bar{\partial}_{i}(p_{\alpha_1\cdots\alpha_k})|^2.
\end{eqnarray*}
From (\ref{19}) it follows that
\begin{eqnarray*}
\sum\limits_{i=1}^3|\bar{\partial}_{i}(p_{\alpha_1\cdots\alpha_k})|^2&=&\sum\limits_{i=1}^3|(\tilde{S}_1^i\cdot\tilde{S}_1+\tilde{S}_2^i\cdot \tilde{S}_2)(p_{\alpha_1\cdots\alpha_k})|^2\\
&=&\sum\limits_{i=1}^3(\tilde{S}_1^i)^2\cdot|\tilde{S}_1(p_{\alpha_1\cdots\alpha_k})|^2+2\sum\limits_{i=1}^3\tilde{S}_1^i\cdot\tilde{S}_2^i\cdot \tilde{S}_1(p_{\alpha_1\cdots\alpha_k})\cdot \tilde{S}_2(p_{\alpha_1\cdots\alpha_k})\\
&&+\sum\limits_{i=1}^3(\tilde{S}_2^i)^2\cdot|\tilde{S}_1(p_{\alpha_1\cdots\alpha_k})|^2.
\end{eqnarray*}
Then the result follows from $\sum\limits_{i=1}^3\tilde{S}^i_k\cdot \tilde{S}^i_l=\delta_{kl}$ for $k,l\in\{1,2\}$.
\endproof
\begin{rem}
Because of Lemma \ref{lem2.7}, now we can say that for any 2-tensor $p$ and $\mathcal{V},\mathcal{W}\in\{\mathcal{T},\mathcal{U},\mathcal{L},\mathcal{S}\}$, the quantity $|\bar{\partial}p|_{\mathcal{V}\mathcal{W}}$ is equivalent to that of (4.5) in \cite{LR2010}. Moreover, $|p|_{\mathcal{V}\mathcal{W}}$ and $|\partial p|_{\mathcal{V}\mathcal{W}}$ are all equivalent to those of (4.3) and (4.4) in \cite{LR2010}.
\end{rem}

\section{Local well-posedness}\label{section3}
From now on, we always assume that $M$ is simply connected. It is easy to check that so is $V_T:=M\times(-T,T)$ with $T\in(0,\infty]$. In such case, the first set of Maxwell equations $d\tilde{F}=0$ is equivalent to $\tilde{F}=dA$ for some 1-form $A$, i.e. $\tilde{F}_{\alpha\beta}=\partial_{\alpha}A_{\beta}-\partial_{\beta}A_{\alpha}$, where $x^0:=t$ and $(x^i,x^0)$ is a natural coordinate chart of $V_T$($(x^i)$ is a local natural coordinate of $M$). Given two Riemannian or pseudo Riemannian metric $g$ and $\hat{e}$ on $V_T$, we suppose that $g$ is in wave gauge with respect to $\hat{e}$. $\hat{D}$ and $\nabla$ are the Levi-Civita connections of $(V_T,\hat{e})$ and $(V_T,g)$ respectively. It is not hard to check that $\tilde{F}_{\alpha\beta}=\hat{D}_{\alpha}A_{\beta}-\hat{D}_{\beta}A_{\alpha}$. For simplicity, we decompose $A$ as
$$A:=A^{\mbox{time}}+A^{\mbox{space}},$$
where
$$
A^{\mbox{time}}(x,t):=A_0(x,t)dt\s\s\s\s\mbox{and}\s\s\s\s A^{\mbox{space}}(x,t):=A_i(x,t)dx^i.
$$

\subsection{Einstein-Maxwell initial data set on simply connected manifold}
In case $M$ is simply connected, the concept of Einstein-Maxwell initial data set should be modified a little. Namely, it is a eight-tuple $(M,\bar{g},K,\bar{A}^{\mbox{space}},\bar{A}_0,\bar{E},\bar{\beta},\bar{N})$, where $\bar{A}^{\mbox{space}}:=\bar{A}_i\cdot dx^i$ is a one form on $M$ and $\bar{A}_0$ is a function on $M$. All the others are the same as that in Subsection \ref{sub2.8}. At this time, $\bar{F}$ in (\ref{3}) should be given by $\bar{F}=d\bar{A}^{\mbox{space}}$.

\subsection{The second set of Maxwell equations in wave and Lorentzian gauges}\label{sub1}
The following computation is obvious
\begin{eqnarray*}
\hat{D}_{\lambda}\hat{D}_{\alpha}A_{\beta}-\hat{D}_{\lambda}\hat{D}_{\beta}A_{\alpha}=\partial_{\lambda}\partial_{\alpha}A_{\beta}-\partial_{\lambda}\partial_{\beta}A_{\alpha}+\hat{\Gamma}^{\theta}_{\lambda\alpha}(\partial_{\beta}A_{\theta}-\partial_{\theta}A_{\beta})+\hat{\Gamma}^{\theta}_{\lambda\beta}(\partial_{\theta}A_{\alpha}-\partial_{\alpha}A_{\theta}),
\end{eqnarray*}
which implies
\begin{eqnarray}\label{7}
\nabla^{\alpha}\tilde{F}_{\alpha\beta}=g^{\alpha\lambda}(\hat{D}_{\lambda}\hat{D}_{\alpha}A_{\beta}-\hat{D}_{\lambda}\hat{D}_{\beta}A_{\alpha})+g^{\alpha\lambda}(\hat{\Gamma}^{\theta}_{\lambda\beta}-\Gamma^{\theta}_{\lambda\beta})\cdot(\hat{D}_{\alpha}A_{\theta}-\hat{D}_{\theta}A_{\alpha}),
\end{eqnarray}
where $\Gamma$ and $\hat{\Gamma}$ are the Christoffel symbols of $g$ and $\hat{e}$ respectively. Here we have used the wave gauge condition.

From Ricci identity it follows that
\begin{eqnarray*}
\hat{D}_{\lambda}\hat{D}_{\beta}A_{\alpha}=\hat{D}_{\beta}\hat{D}_{\lambda}A_{\alpha}+A_{\theta}\cdot\hat{R}^{\theta}_{\alpha\beta\lambda}
\end{eqnarray*}
with
\begin{eqnarray*}
\hat{R}^{\theta}_{\alpha\beta\lambda}:=\partial_{\beta}\hat{\Gamma}^{\theta}_{\lambda\alpha}-\partial_{\lambda}\hat{\Gamma}^{\theta}_{\beta\alpha}+\hat{\Gamma}^{\mu}_{\lambda\alpha}\cdot\hat{\Gamma}^{\theta}_{\beta\mu}-\hat{\Gamma}^{\mu}_{\beta\alpha}\cdot\hat{\Gamma}^{\theta}_{\lambda\mu}.
\end{eqnarray*}
By elementary manipulations (\ref{7}) becomes
\begin{eqnarray}\label{8}
\nabla^{\alpha}\tilde{F}_{\alpha\beta}&=&g^{\alpha\lambda}\hat{D}_{\lambda}\hat{D}_{\alpha}A_{\beta}-g^{\alpha\lambda}\hat{D}_{\beta}\hat{D}_{\lambda}A_{\alpha}-A_{\theta}\cdot\hat{R}^{\theta}_{\alpha\beta\lambda}\cdot g^{\alpha\lambda}\nonumber\\
&&+g^{\alpha\lambda}(\hat{\Gamma}^{\theta}_{\lambda\beta}-\Gamma^{\theta}_{\lambda\beta})\cdot(\hat{D}_{\alpha}A_{\theta}-\hat{D}_{\theta}A_{\alpha}),\\
&=&g^{\alpha\lambda}\hat{D}_{\lambda}\hat{D}_{\alpha}A_{\beta}-\hat{D}_{\beta}(g^{\alpha\lambda}\hat{D}_{\lambda}A_{\alpha})+\hat{D}_{\beta}g^{\alpha\lambda}\cdot\hat{D}_{\lambda}A_{\alpha}-A_{\theta}\cdot\hat{R}^{\theta}_{\alpha\beta\lambda}\cdot g^{\alpha\lambda}\nonumber\\
&&+g^{\alpha\lambda}(\hat{\Gamma}^{\theta}_{\lambda\beta}-\Gamma^{\theta}_{\lambda\beta})\cdot(\hat{D}_{\alpha}A_{\theta}-\hat{D}_{\theta}A_{\alpha}).\nonumber
\end{eqnarray}
We can transform (\ref{8}) into
\begin{eqnarray*}
\nabla^{\alpha}\tilde{F}_{\alpha\beta}=g^{\alpha\lambda}\hat{D}_{\lambda}\hat{D}_{\alpha}A_{\beta}+\hat{D}_{\beta}g^{\alpha\lambda}\cdot\hat{D}_{\lambda}A_{\alpha}-A_{\theta}\cdot\hat{R}^{\theta}_{\alpha\beta\lambda}\cdot g^{\alpha\lambda}+g^{\alpha\lambda}(\hat{\Gamma}^{\theta}_{\lambda\beta}-\Gamma^{\theta}_{\lambda\beta})\cdot(\hat{D}_{\alpha}A_{\theta}-\hat{D}_{\theta}A_{\alpha}),
\end{eqnarray*}
if we assume that $g^{\alpha\lambda}\hat{D}_{\lambda}A_{\alpha}=0$, which is equivalent to $\mbox{div}_gA\equiv g^{\alpha\lambda}\nabla_{\lambda}A_{\alpha}=0$ called the Lorentz gauge condition(the equivalence follows from the wave gauge condition).

It is easy to see
\begin{eqnarray*}
\Gamma^{\lambda}_{\theta\beta}-\hat{\Gamma}^{\lambda}_{\theta\beta}=\frac{1}{2}g^{\lambda\mu}\cdot(\hat{D}_{\theta}g_{\mu\beta}+\hat{D}_{\beta}g_{\theta\mu}-\hat{D}_{\mu}g_{\theta\beta})
\end{eqnarray*}
and
\begin{eqnarray*}
\hat{D}_{\beta}g^{\alpha\lambda}=-g^{\alpha\theta}\cdot\hat{D}_{\beta}g_{\theta\mu}\cdot g^{\mu\lambda}.
\end{eqnarray*}

Therefore, the second set of Maxwell equation in wave and Lorentzian gauges can be written as
\begin{eqnarray}\label{20}
g^{\alpha\lambda}\hat{D}_{\lambda}\hat{D}_{\alpha}A_{\beta}+f_{\beta}(g,A,\hat{D}g,\hat{D}A)=0,
\end{eqnarray}
where
\begin{eqnarray*}
f_{\beta}(g,A,\hat{D}g,\hat{D}A)&:=&-g^{\alpha\theta}\cdot\hat{D}_{\beta}g_{\theta\mu}\cdot g^{\mu\lambda}\cdot\hat{D}_{\lambda}A_{\alpha}-A_{\theta}\cdot\hat{R}^{\theta}_{\alpha\beta\lambda}\cdot g^{\alpha\lambda}\\
&&-\frac{1}{2}g^{\alpha\lambda}g^{\theta\mu}\cdot(\hat{D}_{\lambda}g_{\mu\beta}+\hat{D}_{\beta}g_{\lambda\mu}-\hat{D}_{\mu}g_{\lambda\beta})\cdot(\hat{D}_{\alpha}A_{\theta}-\hat{D}_{\theta}A_{\alpha}).
\end{eqnarray*}
\subsection{Einstein equation in wave gauge}\label{sub2}
Referring to Section 7.4 of Chapter 6 in \cite{CB} we get the coming formula
\begin{eqnarray*}
Ric_{\alpha\beta}=-\frac{1}{2}g^{\lambda\mu}\hat{D}_{\lambda}\hat{D}_{\mu}g_{\alpha\beta}+h_{\alpha\beta}(g,\hat{D}g)+\frac{1}{2}(g_{\alpha\lambda}\hat{D}_{\beta}\hat{F}^{\lambda}+g_{\beta\lambda}\hat{D}_{\alpha}\hat{F}^{\lambda})
\end{eqnarray*}
with
\begin{eqnarray*}
h_{\alpha\beta}(g,\hat{D}g):=P_{\alpha\beta}^{\rho\sigma\gamma\delta\lambda\mu}(g,g^{-1})\cdot\hat{D}_{\rho}g_{\gamma\delta}\cdot\hat{D}_{\sigma}g_{\lambda\mu}+\frac{1}{2}g^{\lambda\mu}\cdot(g_{\alpha\rho}\cdot\hat{R}^{\rho}_{\lambda\beta\mu}+g_{\beta\rho}\cdot\hat{R}^{\rho}_{\lambda\alpha\mu}),
\end{eqnarray*}
where $Ric$ is the Ricci tensor of $g$ and the tensor $P$ is a polynomial in $g$ and $g^{-1}$.

Hence, (\ref{0}) can be reduced to
\begin{eqnarray*}
g^{\lambda\mu}\hat{D}_{\lambda}\hat{D}_{\mu}g_{\alpha\beta}+\tilde{f}_{\alpha\beta}(g,A,\hat{D}g,\hat{D}A)=0,
\end{eqnarray*}
where
\begin{eqnarray*}
\tilde{f}_{\alpha\beta}(g,A,\hat{D}g,\hat{D}A)&:=&2g^{\lambda\mu}\cdot(\hat{D}_{\lambda}A_{\alpha}-\hat{D}_{\alpha}A_{\lambda})\cdot(\hat{D}_{\mu}A_{\beta}-\hat{D}_{\beta}A_{\mu})-2h_{\alpha\beta}(g,\hat{D}g)\\
&&-\frac{1}{n-1}g_{\alpha\beta}g^{\lambda\rho}g^{\mu\sigma}\cdot(\hat{D}_{\lambda}A_{\mu}-\hat{D}_{\mu}A_{\lambda})\cdot(\hat{D}_{\rho}A_{\sigma}-\hat{D}_{\sigma}A_{\rho})\\
&&-g_{\alpha\lambda}\hat{D}_{\beta}\hat{F}^{\lambda}-g_{\beta\lambda}\hat{D}_{\alpha}\hat{F}^{\lambda}.
\end{eqnarray*}
The wave gauge conditions tell us that $\hat{F}=0$. Hence we have
\begin{eqnarray*}
\tilde{f}_{\alpha\beta}(g,A,\hat{D}g,\hat{D}A)&:=&2g^{\lambda\mu}\cdot(\hat{D}_{\lambda}A_{\alpha}-\hat{D}_{\alpha}A_{\lambda})\cdot(\hat{D}_{\mu}A_{\beta}-\hat{D}_{\beta}A_{\mu})-2h_{\alpha\beta}(g,\hat{D}g)\\
&&-\frac{1}{n-1}g_{\alpha\beta}g^{\lambda\rho}g^{\mu\sigma}\cdot(\hat{D}_{\lambda}A_{\mu}-\hat{D}_{\mu}A_{\lambda})\cdot(\hat{D}_{\rho}A_{\sigma}-\hat{D}_{\sigma}A_{\rho}).
\end{eqnarray*}

\subsection{Reducing EM equations to quasi-linear systems on a new bundle over $(V_T,\hat{e})$}\label{sub3.4}
Firstly, we want to construct a new vector bundle $BU$ over $(V_T,\hat{e})$. It is given by
\begin{eqnarray*}
BU:=(T^*V_T\otimes T^*V_T)\times T^*V_T
\end{eqnarray*}
endowed with a metric $\langle\cdot,\cdot\rangle$, where the symbol ``$\times$'' means the Cartesian product of two vector bundles and $T^*V_T$ is the cotangent bundle of $V_T$. More precisely, for any $(g_i,A_i)\in BU$($i=1,2$), we define
\begin{eqnarray*}
\langle(g_1,A_1),(g_2,A_2)\rangle:=((g_1,g_2))+\langle\langle A_1,A_2\rangle\rangle,
\end{eqnarray*}
where
\begin{eqnarray*}
((g_1,g_2)):=(\hat{e})^{\alpha\beta}\cdot(\hat{e})^{\theta\gamma}\cdot(g_1)_{\alpha\theta}\cdot(g_2)_{\beta\gamma}
\end{eqnarray*}
and
\begin{eqnarray*}
\langle\langle A_1,A_2\rangle\rangle:=(\hat{e})^{\alpha\beta}\cdot A_{\alpha}\cdot A_{\beta}
\end{eqnarray*}
with $(\hat{e})^{\alpha\beta}$ being the inverse of $\hat{e}_{\alpha\beta}$. Furthermore, we define a connection $\textbf{\mbox{D}}$ on $BU$ by the following identity
\begin{eqnarray*}
\textbf{\mbox{D}}(g,A):=(\hat{D}g,\hat{D}A)\s\s\s\s\mbox{for any $(g,A)\in BU$.}
\end{eqnarray*}
It is not difficult to check that $\textbf{\mbox{D}}$ is compatible to the metric $\langle\cdot,\cdot\rangle$.

From the discussion in Subsection \ref{sub1} and \ref{sub2} we infer that if $u:=(g,A)\in BU$ satisfies EM equations in wave and Lorentzian gauges, then it is also a solution of the following quasi-linear system
\begin{eqnarray}\label{18}
h^{\lambda\mu}(u)\cdot\textbf{\mbox{D}}_{\lambda}\textbf{\mbox{D}}_{\mu}u+l(u,\textbf{\mbox{D}}u)=0,
\end{eqnarray}
where
\begin{eqnarray*}
h^{\lambda\mu}(u):=g^{\lambda\mu}\s\s\s\s\mbox{and}\s\s\s\s l(u,\textbf{\mbox{D}}u):=\left(\tilde{f}_{\alpha\beta}(u,\textbf{\mbox{D}}u)dx^{\alpha}\otimes dx^{\beta},f_{\theta}(u,\textbf{\mbox{D}}u)dx^{\theta}\right).
\end{eqnarray*}

\subsection{Fixing $\hat{e}$ and determining the initial value of $\textbf{\mbox{D}}_0u$ on $M_0$}\label{sub3.5}
In order to determining $\textbf{\mbox{D}}_0u|_{t=0}$, we must give the values of $(\hat{D}_0g_{ij},\hat{D}_0\beta^i,\hat{D}_0N\equiv\partial_tN,\hat{D}_0A_{\alpha})|_{t=0}$. From now on, we fix $\hat{e}:=dt\otimes dt+e$ with $e$ a Sobolev regular Riemannian metric on $M$. Hence, we get
$$\hat{\Gamma}^{\beta}_{0\alpha}=0,\s\s\s\s\hat{\Gamma}^0_{\alpha\beta}=0\s\s\s\s\mbox{and}\s\s\s\s\hat{\Gamma}^k_{ij}=\tilde{\Gamma}^k_{ij}(\text{$\tilde{\Gamma}$ is the Christoffel symbol of $(M,e)$}),$$
and our goal turns to be determining $(\partial_tg_{ij},\partial_t\beta^i,\partial_tN,\partial_tA_{\alpha})$ on $M_0$. From Chapter 6 of \cite{CB} it follows that they can not be chosen arbitrarily; they should satisfy some restrictions.

By (6.1) of Section 6.1 of Chapter 6 of \cite{CB} we know
\begin{eqnarray}\label{15}
\partial_tg_{ij}|_{t=0}=-2\bar{N}K_{ij}+\bar{g}_{jh}\bar{\nabla}_i\bar{\beta}^h+\bar{g}_{ih}\bar{\nabla}_j\bar{\beta}^h.
\end{eqnarray}

From Lorentzian gauge condition $\mbox{div}_gA=0$ and wave gauge condition we infer that
\begin{eqnarray*}
g^{0\alpha}\partial_tA_{\alpha}=-g^{i\alpha}\hat{D}_iA_{\alpha},
\end{eqnarray*}
implying
\begin{eqnarray*}
-N^{-2}\partial_tA_0+N^{-2}\beta^i\partial_tA_i=-N^{-2}\beta^i\partial_iA_0-(g_*^{ij}-N^{-2}\beta^i\beta^j)\hat{D}_iA_j.
\end{eqnarray*}
Restricting the above identity to $M_0$ yields
\begin{eqnarray}\label{16}
(\partial_tA_0)|_{t=0}=\bar{\beta}^i(\partial_tA_i)|_{t=0}+\bar{\beta}^i\partial_i\bar{A}_0+(\bar{N}^{2}\bar{g}^{ij}-\bar{\beta}^i\bar{\beta}^j)\tilde{D}_i\bar{A}_j,
\end{eqnarray}
where $\bar{A}:=A^{\mbox{space}}|_{t=0}$, $\bar{A}_0:=A^{\mbox{time}}(\partial_t)|_{t=0}$ and $\tilde{D}$ is the Levi-Civita connection of $(M,e)$. Now the problem turns to be how to determine $(\partial_tA_i)|_{t=0}$. Easily, from (\ref{10}) it follows that
\begin{eqnarray*}
(\hat{D}_0A_i-\hat{D}_iA_0)|_{t=0}=\bar{N}\cdot\bar{E}^j\cdot\bar{g}_{ji},
\end{eqnarray*}
which means
\begin{eqnarray}\label{17}
(\partial_tA_i)|_{t=0}=\partial_i\bar{A}_0+\bar{N}\cdot\bar{E}^j\cdot\bar{g}_{ji}.
\end{eqnarray}
Substituting (\ref{17}) into (\ref{16}) yields
\begin{eqnarray}\label{24}
(\partial_tA_0)|_{t=0}=2\bar{\beta}(\bar{A}_0)+\bar{N}\cdot\bar{g}(\bar{E},\bar{\beta})+(\bar{N}^{2}\bar{g}^{ij}-\bar{\beta}^i\bar{\beta}^j)\tilde{D}_i\bar{A}_j.
\end{eqnarray}
In other words, if $\bar{A}_0$ and $\bar{A}$ are given, then $(\partial_tA_0)|_{t=0}$ and $(\partial_tA_i)|_{t=0}$ can be specified via (\ref{16}) and (\ref{17}).

By wave gauge condition $\mbox{tr}_g(\Gamma-\hat{\Gamma})=0$ we obtain
\begin{eqnarray*}
g^{\alpha\beta}\Gamma^0_{\alpha\beta}=0\s\s\s\s\mbox{and}\s\s\s\s g^{\alpha\beta}\Gamma^k_{\alpha\beta}=g^{ij}\tilde{\Gamma}^k_{ij},
\end{eqnarray*}
which are equivalent to
\begin{eqnarray}\label{12}
\partial_tN=\frac{1}{2}Ng_*^{ij}\partial_tg_{ij}-N\mbox{div}_{g_t}\beta+\beta(N)
\end{eqnarray}
and
\begin{eqnarray}\label{13}
\partial_t\beta^k&=&(N^2g_*^{ij}-\beta^i\beta^j)({}^t\Gamma^k_{ij}-\tilde{\Gamma}^k_{ij})-\frac{1}{2}\beta^k\cdot\partial_tg_{ij}\cdot g_*^{ij}+\frac{1}{2}\beta^k\cdot(L_{\beta}g)_{ij}\cdot g_*^{ij}\nonumber\\
&&+N^{-1}\partial_tN\beta^k+\frac{1}{2}g_{*}^{kh}\partial_h\{g(\beta,\beta)-N^2\}-N^{-1}\beta^k\beta(N)\\
&&+\beta^i\cdot\left({}^t\nabla_i\beta^k\right)-\beta^ig_*^{hk}g_{ip}\cdot\left({}^t\nabla_h\beta^p\right).\nonumber
\end{eqnarray}
Here, ${}^t\nabla$ is the Levi-Civita connection of $(M_t,g_t)$ and ${}^t\Gamma^k_{ij}$ is its coefficients, while $L_{\beta}$ is the Lie derivative with respect to $\beta$. Substituting (\ref{12}) into (\ref{13}) yields
\begin{eqnarray}\label{14}
\partial_t\beta^k&=&(N^2g_*^{ij}-\beta^i\beta^j)({}^t\Gamma^k_{ij}-\tilde{\Gamma}^k_{ij})+\frac{1}{2}\beta^k\cdot(L_{\beta}g)_{ij}\cdot g_*^{ij}-\beta^k\mbox{div}_{g_t}\beta\nonumber\\
&&+\frac{1}{2}g_{*}^{kh}\partial_h\{g(\beta,\beta)-N^2\}+\beta^i\cdot\left({}^t\nabla_i\beta^k\right)-\beta^ig_*^{hk}g_{ip}\cdot\left({}^t\nabla_h\beta^p\right).
\end{eqnarray}
Restricting (\ref{12}) and (\ref{13}) to $M_0$ and then substituting (\ref{15}) into them lead to
\begin{eqnarray}\label{25}
\partial_tN|_{t=0}=-\bar{N}^2\cdot\mbox{tr}_{\bar{g}}K+\bar{\beta}(\bar{N})
\end{eqnarray}
and
\begin{eqnarray}\label{26}
\partial_t\beta^k|_{t=0}&=&(\bar{N}^2\bar{g}^{ij}-\bar{\beta}^i\bar{\beta}^j)(\bar{\Gamma}^k_{ij}-\tilde{\Gamma}^k_{ij})+\frac{1}{2}\bar{\beta}^k\cdot\mbox{tr}_{\bar{g}}(L_{\bar{\beta}}\bar{g})-\bar{\beta}^k\mbox{div}_{\bar{g}}\bar{\beta}\\
&&+\frac{1}{2}\bar{g}^{kh}\partial_h\{\bar{g}(\bar{\beta},\bar{\beta})-\bar{N}^2\}+\bar{\beta}^i\cdot\bar{\nabla}_i\bar{\beta}^k-\bar{\beta}^i\bar{g}^{hk}\bar{g}_{ip}\cdot\bar{\nabla}_h\bar{\beta}^p,\nonumber
\end{eqnarray}
where $\bar{\Gamma}^k_{ij}$ is the connection coefficient of $(M,\bar{g})$.
\subsection{Local in time existence and uniqueness in the wave and Lorentzian gauges}
Thanks to Subsection \ref{sub3.4} and \ref{sub3.5}, we have formulated the intrinsic Cauchy problem for Einstein-Maxwell equations in the form of standard PDE analyses. Hence, one can now use the results in Appendix 3 of \cite{CB} to obtain a local in time, global in space, existence and uniqueness theorem in the $\hat{e}$-wave and Lorentzian gauges. The methods we rely on are almost the same as those of Section 7 and 8 of Chapter 6 of \cite{CB}.

\begin{thm}
We take $\hat{e}:=dt\otimes dt+e$ with $e$ a smooth, Sobolev regular metric on a simply connected $n$-dimensional manifold $M$, used to define the Sobolev spaces $H_s$ and $\overset{0}{H_s}$.

Hypotheses on the initial datum sets $(M,\bar{g},K,\bar{A}_0,\bar{A}^{\mbox{space}},\bar{E},\bar{\beta},\bar{N})$ and\\
$(\partial_tg_{ij}|_{t=0}, \partial_tA_0|_{t=0}, \partial_tA_i|_{t=0}, \partial_t\beta^k|_{t=0}, \partial_tN|_{t=0})$:

1. $\bar{g}$ is a Riemannian metric on $M$ uniformly equivalent to $e$ and such that
\begin{eqnarray*}
\tilde{D}\bar{g}\in\overset{0}{H}_{s-1}\s\s\s\s\mbox{and}\s\s\s\s\bar{g}\in\bar{C}^0\s\s\s\s\mbox{with}\s\s\s\s s\in\mathbb{Z}\cap\left(\frac{n}{2}+1,\infty\right),
\end{eqnarray*}
where $\mathbb{Z}$ is the set of all the integers and $\tilde{D}$ is the Levi-Civita connection of $(M,e)$. Furthermore, $(\partial_tg_{ij})|_{t=0}$ is given by (\ref{15}).

2. $K$ is a symmetric 2-tensor on $M$ such that $K\in\overset{0}{H}_{s-1}.$

3. $\bar{A}_0$ belongs to $\bar{C}^0$ and $\partial\bar{A}_0\in\overset{0}{H}_{s-1}$. And $(\partial_tA_0)|_{t=0}$ is given by (\ref{24}).

4. $\bar{A}^{\mbox{space}}\in\bar{C}^0$ and $\tilde{D}(\bar{A}^{\mbox{space}})\in\overset{0}{H}_{s-1}$. Moreover, $(\partial_tA_i)|_{t=0}$ is given by (\ref{17}).

5. $\bar{\beta}\in\bar{C}^0$ and $\tilde{D}\bar{\beta}\in\overset{0}{H}_{s-1}$. Moreover, there exists a positive constant $b$ such that $e(\bar{\beta},\bar{\beta})\leqslant b$. And $\partial_t\beta^k|_{t=0}$ is given by (\ref{26}).

6. $\bar{N}\in\bar{C}^0$ and $\partial\bar{N}\in\overset{0}{H}_{s-1}$. Moreover, there exist two positive constants $C_1$ and $C_2$ such that $C_1\leqslant\bar{N}\leqslant C_2$. Besides, $\partial_tN|_{t=0}$ is given by (\ref{25}).

7. $\bar{E}\in\overset{0}{H}_{s-1}$.

8. $(M,\bar{g},K,\bar{A}_0,\bar{A}^{\mbox{space}},\bar{E},\bar{\beta},\bar{N})$ satisfies the Einstein-Maxwell constraints.

Conclusions:

The initial datum sets admit a development $(V_T,g,A)$ for some $T>0$, such that $A\in\tilde{E}_s(T)$, the spacetime metric $g$ is a regular sliced Lorentzian metric in $\tilde{E}_s(T)$ and $(g,A)$ satisfies on $V_T$ the Einstein-Maxwell equations. Furthermore, $(g,A)$ meets the $\hat{e}$-wave and Lorentzian gauge conditions.

Two such developments in the $\hat{e}$-wave gauge and Lorentzian gauge $(V_T,g_1,A_1)$ and $(V_T,g_2,A_2)$, which are in $\tilde{E}_s(T)$, and which take the same initial values\\
$(\bar{g},K,\bar{A}_0,\bar{A}^{\mbox{space}},\bar{E},\bar{\beta},\bar{N})$ on $M$, coincide on $V_T$.
\end{thm}
\textbf{Sketch of the proof.}
Note that (\ref{18}) are quasi-diagonal, hyperquasi-linear(i.e. $h^{\lambda\mu}$ depends on $u$ but not on $\textbf{\mbox{D}}u$), second-order systems of the type treated in Appendix 3 of \cite{CB}. They satisfies the hypotheses enunciated in that appendix. So the existence and uniqueness theorem for (\ref{18}) then follows.

By Lemma 10.2 of Chapter 6 in \cite{CB} we know that, since the initial datum satisfy the Einstein-Maxwell constraints and
\begin{eqnarray}\label{27}
\hat{F}^{\lambda}|_{t=0}=0,\s\s\s\s (\mbox{div}_gA)|_{t=0}=0,
\end{eqnarray}
the following identities hold true
\begin{eqnarray*}
\partial_t\hat{F}^{\lambda}|_{t=0}=0\s\s\mbox{and}\s\s\partial_t(\mbox{div}_gA)|_{t=0}=0,
\end{eqnarray*}
where it is obvious that the conditions (\ref{17}), (\ref{24}), (\ref{15}), (\ref{25}) and (\ref{26}) lead to (\ref{27}). Furthermore, Lemma 10.1 of Chapter 6 in \cite{CB} tells us that, if $(g,A)$ satisfies (\ref{18}), then $\hat{F}$ and $\mbox{div}_gA$ satisfy a quasi-diagonal system of linear homogeneous differential equations with principal terms the wave equation in the metric $g$. Combining the above two lemmas we arrive at that
$$\hat{F}^{\lambda}=0\s\s\s\s\mbox{and}\s\s\s\s\mbox{div}_gA=0,$$
provided the initial datum satisfy the Einstein-Maxwell constraints and (\ref{27}). Hence, a solution for (\ref{18}), with initial datum satisfying the Einstein-Maxwell constraints and (\ref{27}), is a solution for the full EM system.
\endproof

\section{The EM equations on $\mathbb{R}^4$}\label{section4}
In the sequel, we always assume that $M=\mathbb{R}^3$(i.e. $n=3$) and $e$ is the standard Euclidean metric of $\mathbb{R}^3$. It is well known that in this case, $H_k=\overset{0}{H_k}$ for all the natural number $k$. Given the initial data set $(\mathbb{R}^3;\bar{g},K,\bar{A}^{\mbox{space}},\bar{A}_0,\bar{E},\bar{\beta}\equiv0,\bar{N})$ satisfying the Einstein-Maxwell constraints:
\begin{eqnarray*}
\bar{R}-|K|^2_{\bar{g}}+(\mbox{tr}_{\bar{g}}K)^2=\bar{g}^{ik}\bar{g}^{jl}\bar{F}_{ij}\bar{F}_{kl}/2+3\bar{g}_{ij}\bar{E}^i\bar{E}^j\s\s\s\s\mbox{with}\s\s\s\s\bar{F}:=d\bar{A}^{\mbox{space}},
\end{eqnarray*}

\begin{eqnarray*}
(\mbox{div}_{\bar{g}}K)_i-\partial_i(\mbox{tr}_{\bar{g}}K)=\bar{E}^l\bar{F}_{li},
\end{eqnarray*}

\begin{eqnarray*}
\bar{\nabla}_h\bar{F}_{ij}+\bar{\nabla}_j\bar{F}_{hi}+\bar{\nabla}_i\bar{F}_{jh}=0,
\end{eqnarray*}

and
\begin{eqnarray*}
\bar{\nabla}_i\bar{E}^i=0,
\end{eqnarray*}
we are going to get a solution for the EM equations. Suppose
\begin{eqnarray*}
\bar{g}\in H_{N+1},\s\s K\in H_{N+1},\s\s \bar{E}\in H_{N+1},\s\s \bar{A}^{\mbox{space}}\in H_{N+1},\s\s \bar{A}_0\in H_{N+1},\s\s \bar{N}\in H_{N+1},
\end{eqnarray*}
where the ``$N$" in ``$H_{N+1}$" is the same as that in ``$\mathcal{E}_N$" and is an integer not smaller than 4. Furthermore, we have to assume that $\bar{g}$ is uniformly equivalent to $e$ and $\bar{N}$ is bounded above and below by some positive constants.

In order to satisfy the wave and Lorentzian gauge conditions, we define the initial datum $g_{\mu\nu}|_{t=0}$, $\partial_tg_{\mu\nu}|_{t=0}$, $A_{\alpha}|_{t=0}$, and $\partial_tA_{\alpha}|_{t=0}$ as follows:
\begin{eqnarray}\label{28}
g_{ij}|_{t=0}:=\bar{g}_{ij},\s\s g_{00}|_{t=0}:=-\bar{N}^2,\s\s g_{0i}|_{t=0}:=0,
\end{eqnarray}
\begin{eqnarray}\label{29}
A^{\mbox{space}}|_{t=0}:=\bar{A}^{\mbox{space}},\s\s A^{\mbox{time}}|_{t=0}:=\bar{A}_0dt|_{t=0},
\end{eqnarray}
\begin{eqnarray}\label{30}
\partial_tg_{ij}|_{t=0}:=-2\bar{N}K_{ij},\s\s\partial_tg_{00}|_{t=0}:=2\bar{N}^3\cdot\mbox{tr}_{\bar{g}}K,
\end{eqnarray}
\begin{eqnarray}\label{31}
\partial_tg_{0l}|_{t=0}:=\bar{N}^2\bar{g}^{ij}\partial_j\bar{g}_{il}-\frac{1}{2}\bar{N}^2\bar{g}^{ij}\partial_l\bar{g}_{ij}-\bar{N}\partial_l\bar{N},
\end{eqnarray}
\begin{eqnarray}\label{32}
\partial_tA_0|_{t=0}:=\bar{N}^2\bar{g}^{ij}\partial_i\bar{A}^{\mbox{space}}_j,\s\s\partial_tA_i|_{t=0}:=\partial_i\bar{A}_0+\bar{N}\cdot\bar{E}^j\cdot\bar{g}_{ji}.
\end{eqnarray}
From (\ref{11}) it follows that giving $\partial_tg_{0l}|_{t=0}$ and $\partial_tg_{00}|_{t=0}$ is equivalent to giving $\partial_t\beta^k|_{t=0}$ and $\partial_tN|_{t=0}$.

Now we obtain a solution $(g,A)\in\tilde{E}_{N+1}(T)$ to the EM equations for some $T>0$, which also satisfies the wave and Lorentzian gauge conditions. In this case, (\ref{20}) becomes
\begin{eqnarray*}
\overset{\sim}{\Box}_gA_{\beta}=g^{\alpha\theta}g^{\mu\lambda}\cdot\partial_{\beta}g_{\theta\mu}\cdot\partial_{\lambda}A_{\alpha}-\frac{1}{2}g^{\alpha\lambda}g^{\theta\mu}\cdot(\partial_{\lambda}g_{\mu\beta}+\partial_{\beta}g_{\mu\lambda}-\partial_{\mu}g_{\lambda\beta})\cdot(\partial_{\alpha}A_{\theta}-\partial_{\theta}A_{\alpha})
\end{eqnarray*}
with $\overset{\sim}{\Box}_g:=g^{\alpha\lambda}\partial_{\lambda}\partial_{\alpha}$. And (\ref{0}) turns to be
\begin{eqnarray}\label{21}
Ric_{\mu\nu}&=&g^{\alpha\beta}\cdot(\partial_{\alpha}A_{\mu}-\partial_{\mu}A_{\alpha})\cdot(\partial_{\beta}A_{\nu}-\partial_{\nu}A_{\beta})\nonumber\\
&&-\frac{1}{4}g_{\mu\nu}\cdot\left\{g^{\alpha\rho}g^{\beta\sigma}\cdot(\partial_{\alpha}A_{\beta}-\partial_{\beta}A_{\alpha})\cdot(\partial_{\rho}A_{\sigma}-\partial_{\sigma}A_{\rho})\right\}.
\end{eqnarray}
From (3.17) of \cite{LR2005} it follows that
\begin{eqnarray*}
Ric_{\mu\nu}=-\frac{1}{2}\overset{\sim}{\Box}_gg_{\mu\nu}+\frac{1}{2}\tilde{P}(\partial_\mu g,\partial_{\nu} g)+\frac{1}{2}\tilde{Q}_{\mu\nu}(\partial g,\partial g),
\end{eqnarray*}
where
\begin{eqnarray*}
\tilde{P}(\partial_\mu g,\partial_{\nu} g):=g^{\alpha\alpha'}g^{\beta\beta'}\cdot\left(\frac{1}{4}\partial_{\mu}g_{\beta\beta'}\partial_{\nu}g_{\alpha\alpha'}-\frac{1}{2}\partial_{\nu}g_{\alpha\beta}\partial_{\mu}g_{\alpha'\beta'}\right)
\end{eqnarray*}
and
\begin{eqnarray*}
&&\tilde{Q}_{\mu\nu}(\partial g,\partial g):=g^{\alpha\alpha'}g^{\beta\beta'}\partial_{\alpha}g_{\beta\mu}\partial_{\alpha'}g_{\beta'\nu}-g^{\alpha\alpha'}g^{\beta\beta'}\left(\partial_{\alpha}g_{\beta\mu}\partial_{\beta'}g_{\alpha'\nu}-\partial_{\beta'}g_{\beta\mu}\partial_{\alpha}g_{\alpha'\nu}\right)\\
&&+g^{\alpha\alpha'}g^{\beta\beta'}(\partial_{\mu}g_{\alpha'\beta'}\partial_{\alpha}g_{\beta\nu}-\partial_{\alpha}g_{\alpha'\beta'}\partial_{\mu}g_{\beta\nu})+g^{\alpha\alpha'}g^{\beta\beta'}(\partial_{\nu}g_{\alpha'\beta'}\partial_{\alpha}g_{\beta\mu}-\partial_{\alpha}g_{\alpha'\beta'}\partial_{\nu}g_{\beta\mu})\\
&&+\frac{1}{2}g^{\alpha\alpha'}g^{\beta\beta'}(\partial_{\beta'}g_{\alpha\alpha'}\partial_{\mu}g_{\beta\nu}-\partial_{\mu}g_{\alpha\alpha'}\partial_{\beta'}g_{\beta\nu})+\frac{1}{2}g^{\alpha\alpha'}g^{\beta\beta'}(\partial_{\beta'}g_{\alpha\alpha'}\partial_{\nu}g_{\beta\mu}-\partial_{\nu}g_{\alpha\alpha'}\partial_{\beta'}g_{\beta\mu}).
\end{eqnarray*}
Hence, (\ref{21}) is equivalent to
\begin{eqnarray}\label{22}
\overset{\sim}{\Box}_gg_{\mu\nu}&=&\tilde{P}(\partial_\mu g,\partial_{\nu} g)+\tilde{Q}_{\mu\nu}(\partial g,\partial g)-2g^{\alpha\beta}\cdot(\partial_{\alpha}A_{\mu}-\partial_{\mu}A_{\alpha})\cdot(\partial_{\beta}A_{\nu}-\partial_{\nu}A_{\beta})\nonumber\\
&&+\frac{1}{2}g_{\mu\nu}\cdot\left\{g^{\alpha\rho}g^{\beta\sigma}\cdot(\partial_{\alpha}A_{\beta}-\partial_{\beta}A_{\alpha})\cdot(\partial_{\rho}A_{\sigma}-\partial_{\sigma}A_{\rho})\right\}.
\end{eqnarray}

Define two 2-tensors
\begin{eqnarray*}
h_{\mu\nu}:=g_{\mu\nu}-m_{\mu\nu}\s\s\s\s\mbox{and}\s\s\s\s H^{\mu\nu}:=g^{\mu\nu}-m^{\mu\nu},
\end{eqnarray*}
where $m^{\mu\nu}$ and $g^{\mu\nu}$ are the inverses of $m_{\mu\nu}$ and $g_{\mu\nu}$ respectively(Recall that $m_{\mu\nu}$ is the Minkowski metric of $\mathbb{R}^4$). We want to obtain the equation of $h_{\mu\nu}$. To this end, the following lemma is necessary.
\begin{lem}\label{lem4.1}
For small $h$, we have
\begin{eqnarray*}
H^{\mu\nu}=-h^{\mu\nu}+O^{\mu\nu}(h^2),
\end{eqnarray*}
where $h^{\mu\nu}:=m^{\mu\mu'}m^{\nu\nu'}h_{\mu'\nu'}$ and $O^{\mu\nu}(h^2)$ means a 2-tensor vanishing to the second order at $h=0$.
\end{lem}
\textbf{Proof.} Hereafter, we always raise and pull down indices via $m^{\mu\nu}$ and $m_{\mu\nu}$. Since $g_{\mu\nu}=m_{\mu\nu}+h_{\mu\nu}$, regarding it as a matrix yields
\begin{eqnarray*}
g_{\mu\nu}=m_{\mu\alpha}\cdot(\delta^{\alpha}_{\nu}+h^{\alpha}_{\nu}),
\end{eqnarray*}
which implies
\begin{eqnarray*}
(g_{\mu\nu})^{-1}=(\delta^{\alpha}_{\nu}+h^{\alpha}_{\nu})^{-1}\cdot(m_{\mu\alpha})^{-1}.
\end{eqnarray*}
This means
\begin{eqnarray*}
g^{\mu\nu}=(\delta^{\alpha}_{\nu}+h^{\alpha}_{\nu})^{-1}\cdot m^{\mu\alpha}.
\end{eqnarray*}
It is easy to check that, if $h$ is small,
\begin{eqnarray*}
(\delta^{\alpha}_{\nu}+h^{\alpha}_{\nu})^{-1}=\delta^{\nu}_{\alpha}-h^{\nu}_{\alpha}+\sum\limits_{i=0}^{\infty}(-1)^i\cdot h^{\nu}_{\kappa_0}\cdot h^{\kappa_0}_{\kappa_1}\cdots h^{\kappa_{i-1}}_{\kappa_i}\cdot h^{\kappa_i}_{\alpha}.
\end{eqnarray*}
Denoting
$$O^{\mu\nu}(h^2):=m^{\mu\alpha}\cdot\sum\limits_{i=0}^{\infty}(-1)^i\cdot h^{\nu}_{\kappa_0}\cdot h^{\kappa_0}_{\kappa_1}\cdots h^{\kappa_{i-1}}_{\kappa_i}\cdot h^{\kappa_i}_{\alpha}$$
leads to the result of this lemma.
\endproof

\begin{rem}
By similar way, we can also get that for small $H$,
\begin{eqnarray*}
h_{\mu\nu}=-H_{\mu\nu}+O_{\mu\nu}(H^2),
\end{eqnarray*}
where $H_{\mu\nu}:=H^{\alpha\beta}m_{\alpha\mu}m_{\beta\nu}$ and $O_{\mu\nu}(H^2)$ means a 2-tensor vanishing to the second order at $H=0$.
\end{rem}

By Lemma 3.2 of \cite{LR2005} we know that if $h$ is small, (\ref{22}) is equivalent to
\begin{eqnarray}\label{23}
\overset{\sim}{\Box}_gh_{\mu\nu}&=&P(\partial_{\mu}h,\partial_{\nu}h)+Q_{\mu\nu}(\partial h,\partial h)+G_{\mu\nu}(h)(\partial h,\partial h)\nonumber\\
&&+2(h^{\alpha\beta}-m^{\alpha\beta})(\partial_{\alpha}A_{\mu}-\partial_{\mu}A_{\alpha})(\partial_{\beta}A_{\nu}-\partial_{\nu}A_{\beta})\\
&&-2O^{\alpha\beta}(h^2)\cdot(\partial_{\alpha}A_{\mu}-\partial_{\mu}A_{\alpha})(\partial_{\beta}A_{\nu}-\partial_{\nu}A_{\beta})\nonumber\\
&&+\frac{1}{2}\{m^{\alpha\rho}m^{\beta\sigma}h_{\mu\nu}+m^{\alpha\rho}m^{\beta\sigma}m_{\mu\nu}-h^{\alpha\rho}m^{\beta\sigma}m_{\mu\nu}-m^{\alpha\rho}h^{\beta\sigma}m_{\mu\nu}\nonumber\\
&&+O^{\alpha\rho\beta\sigma}_{\mu\nu}(h^2)\}(\partial_{\alpha}A_{\beta}-\partial_{\beta}A_{\alpha})(\partial_{\rho}A_{\sigma}-\partial_{\sigma}A_{\rho}),\nonumber\\
&:=&F_{\mu\nu}\nonumber
\end{eqnarray}
where $O^{\alpha\rho\beta\sigma}_{\mu\nu}(h^2)$ vanishes to the second order at $h=0$,
\begin{eqnarray*}
P(\partial_{\mu}h,\partial_{\nu}h):=m^{\alpha\alpha'}m^{\beta\beta'}\cdot\left(\frac{1}{4}\partial_{\mu}h_{\beta\beta'}\partial_{\nu}h_{\alpha\alpha'}-\frac{1}{2}\partial_{\nu}h_{\alpha\beta}\partial_{\mu}h_{\alpha'\beta'}\right)
\end{eqnarray*}
and
\begin{eqnarray*}
&&Q_{\mu\nu}(\partial h,\partial h):=m^{\alpha\alpha'}m^{\beta\beta'}\partial_{\alpha}h_{\beta\mu}\partial_{\alpha'}h_{\beta'\nu}-m^{\alpha\alpha'}m^{\beta\beta'}\left(\partial_{\alpha}h_{\beta\mu}\partial_{\beta'}h_{\alpha'\nu}-\partial_{\beta'}h_{\beta\mu}\partial_{\alpha}h_{\alpha'\nu}\right)\\
&&+m^{\alpha\alpha'}m^{\beta\beta'}(\partial_{\mu}h_{\alpha'\beta'}\partial_{\alpha}h_{\beta\nu}-\partial_{\alpha}h_{\alpha'\beta'}\partial_{\mu}h_{\beta\nu})+m^{\alpha\alpha'}m^{\beta\beta'}(\partial_{\nu}h_{\alpha'\beta'}\partial_{\alpha}h_{\beta\mu}-\partial_{\alpha}h_{\alpha'\beta'}\partial_{\nu}h_{\beta\mu})\\
&&+\frac{1}{2}m^{\alpha\alpha'}m^{\beta\beta'}(\partial_{\beta'}h_{\alpha\alpha'}\partial_{\mu}h_{\beta\nu}-\partial_{\mu}h_{\alpha\alpha'}\partial_{\beta'}h_{\beta\nu})+\frac{1}{2}m^{\alpha\alpha'}m^{\beta\beta'}(\partial_{\beta'}h_{\alpha\alpha'}\partial_{\nu}h_{\beta\mu}-\partial_{\nu}h_{\alpha\alpha'}\partial_{\beta'}h_{\beta\mu}),
\end{eqnarray*}
is a null form and $G_{\mu\nu}(h)(\partial h,\partial h)$ is a quadratic form in $\partial h$ with coefficients smoothly dependent on $h$ and vanishing when $h$ vanishes, i.e. $G_{\mu\nu}(0)(\partial h,\partial h)=0$.

Using Lemma \ref{lem4.1} again we get
\begin{eqnarray}\label{36}
\overset{\sim}{\Box}_gA_{\beta}&=&\{-m^{\alpha\theta}h^{\mu\lambda}-h^{\alpha\theta}m^{\mu\lambda}+m^{\alpha\theta}m^{\mu\lambda}+O^{\alpha\theta\mu\lambda}(h^2)\}\partial_{\beta}h_{\theta\mu}\partial_{\lambda}A_{\alpha}\\
&&+\frac{1}{2}\{-m^{\alpha\lambda}h^{\theta\mu}-h^{\alpha\lambda}m^{\theta\mu}+m^{\alpha\lambda}m^{\theta\mu}+O^{\alpha\lambda\theta\mu}(h^2)\}\nonumber\\
&&\times(\partial_{\lambda}h_{\mu\beta}+\partial_{\beta}h_{\mu\lambda}-\partial_{\mu}h_{\lambda\beta})(\partial_{\alpha}A_{\theta}-\partial_{\theta}A_{\alpha})\nonumber\\
&:=&J_{\beta},\nonumber
\end{eqnarray}
provided $h$ is sufficiently small.

\section{Beginning of the proof of Theorem \ref{thm6.1}}\label{section5}
As described in the introduction, $T$ is the maximal existence time of the solution $(g,A)$ and assumed to be finite. We have defined the time $T^*$ to be
\begin{eqnarray}\label{38}
T^*:=\sup\{T_0:\exists C=C(T),\,s.t.\s\forall t\in[0,T_0),\s\mathcal{E}_N(t)\leqslant2C\varepsilon\},
\end{eqnarray}
where $\mathcal{E}_N$ is given by (\ref{95}). Our goal is to show that if $\varepsilon>0$ is small enough, then the inequality in (\ref{38}) implies the same inequality with $2C$ replaced by $C$ for all $t<T^*$.

The first step is to derive the preliminary decay estimates for $h$ and $A$ under the assumption (\ref{38}). In order to make readers easy to understand the proof, we first let $t$ appear in the inequalities below.
\begin{thm}\label{thm6.2}
Let $h$ and $A$ verify the inequality in (\ref{38}). Then
\begin{equation}\label{39}
|\partial Z^Ih(x,t)|+|\partial Z^IA(x,t)|\lesssim\left\{
\begin{array}{rcl}
\varepsilon(1+t+|q|)^{-1}(1+|q|)^{-1-\gamma}&        &\text{$q>0$,}\\
\varepsilon(1+t+|q|)^{-1}(1+|q|)^{-1/2}&         &\text{$q<0$,}
\end{array} \right.
|I|\leqslant N-3.
\end{equation}
Furthermore
\begin{equation}\label{40}
|Z^Ih(x,t)|+|Z^IA(x,t)|\lesssim
\left\{
\begin{array}{rcl}
\varepsilon(1+t+|q|)^{-1}(1+|q|)^{-\gamma}&        &\text{$q>0$,}\\
\varepsilon(1+t+|q|)^{-1}(1+|q|)^{1/2}&         &\text{$q<0$,}
\end{array} \right.
|I|\leqslant N-3.
\end{equation}
And
\begin{equation}\label{41}
|\bar{\partial} Z^Ih(x,t)|+|\bar{\partial} Z^IA(x,t)|\lesssim
\left\{
\begin{array}{rcl}
\varepsilon(1+t+|q|)^{-2}(1+|q|)^{-\gamma}&        &\text{$q>0$,}\\
\varepsilon(1+t+|q|)^{-2}(1+|q|)^{1/2}&         &\text{$q<0$,}
\end{array} \right.
|I|\leqslant N-4.
\end{equation}
\end{thm}
\textbf{Proof.} Estimate (\ref{39}) follows from the weighted Sobolev inequality of Proposition 14.1 of \cite{LR2010}. In particular,
\begin{eqnarray*}
|\partial Z^Ih(x,0)|+|\partial Z^IA(x,0)|\lesssim\varepsilon\cdot(1+|x|)^{-2-\gamma}\s\s\mbox{for all $|I|\leqslant N-3$.}
\end{eqnarray*}
Of course, we have the following
\begin{eqnarray*}
|\partial Z^Ih(x,0)|+|\partial Z^IA(x,0)|\lesssim\varepsilon\cdot(1+|x|)^{-2-\gamma}\s\s\mbox{for all $|I|\leqslant N-4$.}
\end{eqnarray*}
Since $|Z\phi|\lesssim(1+t+|x|)|\partial\phi|$ for any smooth function $\phi$, we say
\begin{eqnarray*}
|Z^Ih(x,0)|+|Z^IA(x,0)|\lesssim\varepsilon\cdot(1+|x|)^{-1-\gamma}\s\s\mbox{for all $|I|\leqslant N-3$.}
\end{eqnarray*}
The estimate (\ref{40}) for $q>0$ follows by integrating (\ref{39}) from the hyperplane $t=0$ along the lines with $t+r$ and $\omega:=x/|x|$ fixed:
\begin{eqnarray*}
|Z^IA(r\omega,t)|+|Z^Ih(r\omega,t)|&\leqslant&\int_r^{t+r}\{|\partial Z^Ih(\rho\omega,t+r-\rho)|+|\partial Z^IA(\rho\omega,t+r-\rho)|\}d\rho\\
&&+|Z^Ih((t+r)\omega,0)|+|Z^IA((t+r)\omega,0)|\\
&\lesssim&\int_r^{t+r}\frac{\varepsilon d\rho}{(1+t+r)(1+|2\rho-t-r|)^{1+\gamma}}+\frac{\varepsilon}{(1+t+r)^{1+\gamma}}\\
&\lesssim&\frac{\varepsilon}{1+t+r}\int_r^{\infty}\frac{d\rho}{(1+|\rho-t|)^{1+\gamma}}+\frac{\varepsilon}{(1+t+|q|)^{1+\gamma}}\\
&\lesssim&\varepsilon(1+t+|q|)^{-1}(1+|q|)^{-\gamma}.
\end{eqnarray*}
A similar argument yields (\ref{40}) for $q<0$. Inequality (\ref{41}) follows from (\ref{40}) using (5.5) of \cite{LR2010}.
\endproof

Recall that in this article, the constant $C$ depends upon $T$. Hence, we give the following theorem which is equivalent to Theorem \ref{thm6.2}.
\begin{thm}\label{thm6.2'}
Let $h$ and $A$ verify the inequality in (\ref{38}). Then
\begin{equation}\label{39'}
|\partial Z^Ih(x,t)|+|\partial Z^IA(x,t)|\lesssim\left\{
\begin{array}{rcl}
\varepsilon(1+|q|)^{-2-\gamma}&        &\text{$q>0$,}\\
\varepsilon(1+|q|)^{-3/2}&         &\text{$q<0$,}
\end{array} \right.
|I|\leqslant N-3.
\end{equation}
Furthermore
\begin{equation}\label{40'}
|Z^Ih(x,t)|+|Z^IA(x,t)|\lesssim
\left\{
\begin{array}{rcl}
\varepsilon(1+|q|)^{-1-\gamma}&        &\text{$q>0$,}\\
\varepsilon(1+|q|)^{-1/2}&         &\text{$q<0$,}
\end{array} \right.
|I|\leqslant N-3.
\end{equation}
And
\begin{equation}\label{41'}
|\bar{\partial} Z^Ih(x,t)|+|\bar{\partial} Z^IA(x,t)|\lesssim
\left\{
\begin{array}{rcl}
\varepsilon(1+|q|)^{-2-\gamma}&        &\text{$q>0$,}\\
\varepsilon(1+|q|)^{-3/2}&         &\text{$q<0$,}
\end{array} \right.
|I|\leqslant N-4.
\end{equation}
\end{thm}
\subsection{Estimates for the inhomogeneous terms $F_{\mu\nu}$ and $J_{\beta}$}
(\ref{40'}) tells us that $|Z^Ih|+|Z^IA|\leqslant1/2$, provided $\varepsilon$ is small enough and $|I|\leqslant N-3$. The smallness play a key role in the sequel.
\begin{pro}\label{pro6.3}
Assume that $h=g-m$ and $A$ satisfy the inequality in (\ref{38}). Let $F_{\mu\nu}$ and $J_{\beta}$ be as in (\ref{23}) and (\ref{36}) respectively. Then we have
\begin{eqnarray}\label{47}
|Z^IF|&\lesssim&\sum\limits_{|J|+|K|\leqslant|I|}(|\partial Z^Jh|_{\mathcal{T}\mathcal{U}}\cdot|\partial Z^Kh|_{\mathcal{T}\mathcal{U}}+|\bar{\partial}Z^Jh|\cdot|\partial Z^Kh|)\nonumber\\
&&+\sum\limits_{|J|+|K|\leqslant|I|-1}|\partial Z^Jh|_{\mathcal{L}\mathcal{T}}\cdot|\partial Z^Kh|+\sum\limits_{|J|+|K|\leqslant|I|-2}|\partial Z^Jh|\cdot|\partial Z^Kh|\\
&&+\sum\limits_{|J_1|+|J_2|+|J_3|\leqslant|I|}|Z^{J_3}h|\cdot|\partial Z^{J_2}h|\cdot|\partial Z^{J_1}h|+\sum\limits_{|J_1|+|J_2|\leqslant|I|}|\partial Z^{J_1}A|\cdot|\partial Z^{J_2}A|\nonumber
\end{eqnarray}
and
\begin{eqnarray}\label{48}
|Z^IJ|\lesssim\sum\limits_{|I_1|+|I_2|\leqslant|I|}|\partial Z^{I_1}h|\cdot|\partial Z^{I_2}A|,
\end{eqnarray}
provided $|Z^Jh|\leqslant\tilde{C}<1$ for all multi-indices $|J|\leqslant|I|$ and vector fields $Z\in\mathfrak{L}$. Here the ``$J$" in (\ref{48}) is the same as that in (\ref{36}).
\end{pro}
\textbf{Proof.} For simplicity, we only show (\ref{47}). Another estimate follows from the same approach.

Reviewing the definition of $F$ gives
\begin{eqnarray*}
Z^IF_{\mu\nu}&=&Z^I\{P(\partial_{\mu}h,\partial_{\nu}h)+Q_{\mu\nu}(\partial h,\partial h)+G_{\mu\nu}(h)(\partial h,\partial h)\}\\
&&+\mbox{Term}_1+\mbox{Term}_2+\mbox{Term}_3+\mbox{Term}_4,
\end{eqnarray*}
where
\begin{eqnarray*}
\mbox{Term}_1:=2Z^I\{(h^{\alpha\beta}-m^{\alpha\beta})(\partial_{\alpha}A_{\mu}-\partial_{\mu}A_{\alpha})(\partial_{\beta}A_{\nu}-\partial_{\nu}A_{\beta})\},
\end{eqnarray*}
\begin{eqnarray*}
\mbox{Term}_2:=-2Z^I\{O^{\alpha\beta}(h^2)\cdot(\partial_{\alpha}A_{\mu}-\partial_{\mu}A_{\alpha})(\partial_{\beta}A_{\nu}-\partial_{\nu}A_{\beta})\},
\end{eqnarray*}
\begin{eqnarray*}
\mbox{Term}_3&:=&\frac{1}{2}Z^I\big\{(m^{\alpha\rho}m^{\beta\sigma}h_{\mu\nu}+m^{\alpha\rho}m^{\beta\sigma}m_{\mu\nu}\nonumber\\
&&-h^{\alpha\rho}m^{\beta\sigma}m_{\mu\nu}-m^{\alpha\rho}h^{\beta\sigma}m_{\mu\nu})(\partial_{\alpha}A_{\beta}-\partial_{\beta}A_{\alpha})(\partial_{\rho}A_{\sigma}-\partial_{\sigma}A_{\rho})\big\},\nonumber
\end{eqnarray*}
and
\begin{eqnarray*}
\mbox{Term}_4:=Z^I\left\{O^{\alpha\rho\beta\sigma}_{\mu\nu}(h^2)(\partial_{\alpha}A_{\beta}-\partial_{\beta}A_{\alpha})(\partial_{\rho}A_{\sigma}-\partial_{\sigma}A_{\rho})\right\}.
\end{eqnarray*}

(9.28) of \cite{LR2010} tells us that, if $|Z^Jh|\leqslant\tilde{C}<1$ for all multi-indices $|J|\leqslant|I|$ and vector fields $Z\in\mathfrak{L}$, one will obtain
\begin{eqnarray}\label{45}
|Z^I(P+Q+G)|&\lesssim&\sum\limits_{|J|+|K|\leqslant|I|}(|\partial Z^Jh|_{\mathcal{T}\mathcal{U}}\cdot|\partial Z^Kh|_{\mathcal{T}\mathcal{U}}+|\bar{\partial}Z^Jh|\cdot|\partial Z^Kh|)\nonumber\\
&&+\sum\limits_{|J|+|K|\leqslant|I|-1}|\partial Z^Jh|_{\mathcal{L}\mathcal{T}}\cdot|\partial Z^Kh|+\sum\limits_{|J|+|K|\leqslant|I|-2}|\partial Z^Jh|\cdot|\partial Z^Kh|\\
&&+\sum\limits_{|J_1|+|J_2|+|J_3|\leqslant|I|}|Z^{J_3}h|\cdot|\partial Z^{J_2}h|\cdot|\partial Z^{J_1}h|.\nonumber
\end{eqnarray}
Moreover, it is easy to get
\begin{eqnarray*}
|\mbox{Term}_1|\lesssim\sum\limits_{|J_1|+|J_2|\leqslant|I|}|Z^{J_1}\partial A|\cdot|Z^{J_2}\partial A|.
\end{eqnarray*}
From induction argument it follows that for any multi-index $I$, there exist a set of universal constants $\{C_J:0\leqslant|J|\leqslant|I|\}$ such that
\begin{eqnarray}\label{46}
Z^I\partial_{\alpha}=\sum\limits_{0\leqslant|J|\leqslant|I|}C_J\cdot\partial_{\alpha}Z^J,
\end{eqnarray}
implying
\begin{eqnarray*}
|\mbox{Term}_1|\lesssim\sum\limits_{|J_1|+|J_2|\leqslant|I|}|\partial Z^{J_1}A|\cdot|\partial Z^{J_2}A|.
\end{eqnarray*}
The same methods leads to
\begin{eqnarray*}
|\mbox{Term}_2|\lesssim\sum\limits_{|J_1|+|J_2|\leqslant|I|}|\partial Z^{J_1}A|\cdot|\partial Z^{J_2}A|,
\end{eqnarray*}
\begin{eqnarray*}
|\mbox{Term}_3|\lesssim\sum\limits_{|J_1|+|J_2|\leqslant|I|}|\partial Z^{J_1}A|\cdot|\partial Z^{J_2}A|
\end{eqnarray*}
and
\begin{eqnarray*}
|\mbox{Term}_4|\lesssim\sum\limits_{|J_1|+|J_2|\leqslant|I|}|\partial Z^{J_1}A|\cdot|\partial Z^{J_2}A|.
\end{eqnarray*}
Then the result of this proposition follows.
\endproof

\section{Energy estimates for EM equations}\label{section6}
In this section we prove the following result.

\begin{thm}\label{thm8.1}
Let $h_{\mu\nu}=g_{\mu\nu}-m_{\mu\nu}$ and $A_{\beta}$ be a local in time solution to (\ref{23}) and (\ref{36}) respectively satisfying the wave and Lorentzian gauge conditions on the interval $[0,T^*)$. Suppose also that we have fixed some constants $\gamma\in(0,1/2)$ and $\mu\in(0,1/2)$. Assume that we have the following estimates for $t\in[0,T^*)$ and all multi-indices $|I|\leqslant N-4$:
\begin{eqnarray}\label{75}
|\partial H|_{\mathcal{T}\mathcal{U}}+\frac{|H|_{\mathcal{L}\mathcal{T}}}{1+|q|}+\frac{|ZH|_{\mathcal{L}\mathcal{L}}}{1+|q|}\lesssim\varepsilon,
\end{eqnarray}
\begin{eqnarray}\label{76}
|\partial Z^Ih|+\frac{|Z^Ih|}{1+|q|}+|\bar{\partial}Z^Ih|\leqslant\left\{
\begin{array}{rcl}
&C\varepsilon(1+|q|)^{-2-\gamma}        &\quad\text{when $q>0$,}\\
&C\varepsilon(1+|q|)^{-3/2}         &\quad\text{when $q<0$,}
\end{array} \right.
\end{eqnarray}
\begin{eqnarray}\label{78}
|\partial Z^IA|+\frac{|Z^IA|}{1+|q|}+|\bar{\partial}Z^IA|\leqslant\left\{
\begin{array}{rcl}
&C\varepsilon(1+|q|)^{-2-\gamma}        &\quad\text{when $q>0$,}\\
&C\varepsilon(1+|q|)^{-3/2}         &\quad\text{when $q<0$}
\end{array} \right.
\end{eqnarray}
and
\begin{eqnarray}\label{77}
\mathcal{E}_N(0)\leqslant\varepsilon.
\end{eqnarray}
Then there is a positive constant $C'$ dependent of $T$ such that we have the energy estimate
\begin{eqnarray*}
\mathcal{E}_N(t)\leqslant C'\varepsilon^2,
\end{eqnarray*}
for all $t\in[0,T^*)$.
\end{thm}

Assuming the conclusions of Theorem \ref{thm8.1} for a moment we finish the proof of the main Theorem \ref{thm6.1}.
\subsection{End of the proof of Theorem \ref{thm6.1}}
Recall that $T^*$ was defined as the maximal time with the property that the bound
$$\mathcal{E}_N(t)\leqslant2C\varepsilon$$
holds for all $t\in[0,T^*)$. Direct check shows that the estimates of Theorem \ref{thm6.2'} imply the assumption (\ref{75})-(\ref{78}). The conclusion of Theorem \ref{thm8.1} states that the energy
$$\mathcal{E}_N(t)\leqslant C'\varepsilon^2,\s\s\s\s\forall t\in[0,T^*).$$
Thus choosing a sufficiently small $\varepsilon>0$ we can show that $\mathcal{E}_N(t)\leqslant C\varepsilon$ thus contracting the maximality of $T^*$ and consequently proving that $(g,A)$ is a global solution. Therefore, it remains to prove Theorem \ref{thm8.1}.
\subsection{Proof of Theorem \ref{thm8.1}}
Recall that $h_{\mu\nu}$ and $A_{\beta}$ satisfy the wave equations $\overset{\sim}{\Box}_gh_{\mu\nu}=F_{\mu\nu}$ and $\overset{\sim}{\Box}_gA_{\beta}=J_{\beta}$ respectively. Our goal is to compute the energy norms of $Z^Ih$ and $Z^IA$, where $Z\in\mathfrak{L}$.

From (11.10) of \cite{LR2010} it follows that
$$
\overset{\sim}{\Box}_gZ^Ih_{\mu\nu}=F^{I}_{\mu\nu},
$$
where
$$
F^{I}:=\hat{Z}^IF-D^I,\s\s\s\s D^I:=\hat{Z}^I\overset{\sim}{\Box}_gh^1-\overset{\sim}{\Box}_gZ^Ih^1,\s\s\s\s \hat{Z}:=Z+c_Z.
$$
Similarly, we also have
$$\overset{\sim}{\Box}_gZ^IA=\hat{Z}^IJ-W^I,$$
where $W^I:=\hat{Z}^I\overset{\sim}{\Box}_gA-\overset{\sim}{\Box}_gZ^IA$. (11.13) of \cite{LR2010} tells us
\begin{eqnarray}\label{80}
&&\int_{\Sigma_t}|\partial Z^Ih|^2w+\int_0^t\int_{\Sigma_{\tau}}|\bar{\partial}Z^Ih|^2w'\\
&\lesssim&\int_{\Sigma_0}|\partial Z^I h|^2w+\int_0^t\int_{\Sigma_{\tau}}\Big\{\frac{\varepsilon|\partial Z^Ih|^2}{1+\tau}w+\frac{w\cdot(1+\tau)}{\varepsilon}(|\hat{Z}^IF|^2+|D^I|^2)\Big\},\nonumber\\
&\lesssim&\int_{\Sigma_0}|\partial Z^I h|^2w+\int_0^t\int_{\Sigma_{\tau}}\Big\{\varepsilon|\partial Z^Ih|^2w+\frac{w}{\varepsilon}(|\hat{Z}^IF|^2+|D^I|^2)\Big\},\nonumber
\end{eqnarray}
where $\Sigma_t:=\mathbb{R}^3\times\{t\}$. Applying the same methods yields
\begin{eqnarray}\label{81}
&&\int_{\Sigma_t}|\partial Z^IA|^2w+\int_0^t\int_{\Sigma_{\tau}}|\bar{\partial}Z^IA|^2w'\\
&\lesssim&\int_{\Sigma_0}|\partial Z^IA|^2w+\int_0^t\int_{\Sigma_{\tau}}\left\{\varepsilon|\partial Z^IA|^2w+\frac{w}{\varepsilon}(|\hat{Z}^IJ|^2+|W^I|^2)\right\}.\nonumber
\end{eqnarray}
We begin with the following estimates on the inhomogeneous terms $F$ and $J$.
\begin{lem}\label{lem8.2}
Under the assumptions of Theorem \ref{thm8.1}, we have
\begin{eqnarray}\label{82}
|Z^IF|&\lesssim&\sum\limits_{|K|\leqslant|I|}\left\{\varepsilon|\partial Z^Kh|+\varepsilon(1+|q|)^{-3/2}|\bar{\partial}Z^Kh|+\varepsilon^2\frac{|Z^Kh|}{(1+|q|)^3}\right\}\nonumber\\
&&+\varepsilon(1+|q|)^{-1}\sum\limits_{|K|\leqslant|I|}|\partial Z^KA|
\end{eqnarray}
and
\begin{equation}\label{86}
|Z^IJ|\lesssim\sum\limits_{|K|\leqslant|I|}|\partial Z^Kh|\cdot\varepsilon(1+|q|)^{-1}.
\end{equation}
\end{lem}
\textbf{Proof.} According to Proposition \ref{pro6.3}
\begin{eqnarray*}
|Z^IF|\lesssim\mbox{Term}+\sum\limits_{|J_1|+|J_2|\leqslant|I|}|\partial Z^{J_1}A|\cdot|\partial Z^{J_2}A|,
\end{eqnarray*}
where
\begin{eqnarray*}
\mbox{Term}&:=&\sum\limits_{|J|+|K|\leqslant|I|}(|\partial Z^Jh|_{\mathcal{T}\mathcal{U}}\cdot|\partial Z^Kh|_{\mathcal{T}\mathcal{U}}+|\bar{\partial}Z^Jh|\cdot|\partial Z^Kh|)\nonumber\\
&&+\sum\limits_{|J|+|K|\leqslant|I|-1}|\partial Z^Jh|_{\mathcal{L}\mathcal{T}}\cdot|\partial Z^Kh|+\sum\limits_{|J|+|K|\leqslant|I|-2}|\partial Z^Jh|\cdot|\partial Z^Kh|\\
&&+\sum\limits_{|J_1|+|J_2|+|J_3|\leqslant|I|}|Z^{J_3}h|\cdot|\partial Z^{J_2}h|\cdot|\partial Z^{J_1}h|,
\end{eqnarray*}
and
\begin{eqnarray*}
|Z^IJ|\lesssim\sum\limits_{|I_1|+|I_2|\leqslant|I|}|\partial Z^{I_1}h|\cdot|\partial Z^{I_2}A|.
\end{eqnarray*}
From (\ref{76}) it follows that
\begin{equation}\label{84}
\mbox{Term}\lesssim\sum\limits_{|K|\leqslant|I|}\left\{\varepsilon|\partial Z^Kh|+\varepsilon(1+|q|)^{-3/2}|\bar{\partial}Z^Kh|+\varepsilon^2\frac{|Z^Kh|}{(1+|q|)^3}\right\}.
\end{equation}
In addition, (\ref{78}) implies
\begin{eqnarray}\label{87}
|\partial Z^KA|\lesssim\varepsilon(1+|q|)^{-1}\s\s\s\s\forall|k|\leqslant|I|.
\end{eqnarray}
This means
\begin{eqnarray}\label{83}
\sum\limits_{|I_1|+|I_2|\leqslant|I|}|\partial Z^{I_1}A|\cdot|\partial Z^{I_2}A|\lesssim\varepsilon(1+|q|)^{-1}\sum\limits_{|K|\leqslant|I|}|\partial Z^KA|.
\end{eqnarray}
Combining (\ref{84}) with (\ref{83}) yields (\ref{82}).

Note that
\begin{eqnarray}\label{85}
|Z^IJ|\lesssim\sum\limits_{|I_1|+|I_2|\leqslant|I|}|\partial Z^{I_1}h|\cdot|\partial Z^{I_2}A|.
\end{eqnarray}
Hence (\ref{86}) follows from (\ref{87}) immediately.
\endproof
\begin{lem}\label{lem8.3}
Under the assumptions of Theorem \ref{thm8.1}, one can get
\begin{eqnarray*}
&&\varepsilon^{-1}\int_0^t\int|Z^IF|^2wdxd\tau\\
&\lesssim&\sum\limits_{|K|\leqslant|I|}\varepsilon\int_0^t\int\left\{|\partial Z^Kh|^2w+|\bar{\partial} Z^Kh|^2w'\right\}dxd\tau+\varepsilon\sum\limits_{|K|\leqslant|I|}\int_0^t\int|\partial Z^KA|^2wdxd\tau,
\end{eqnarray*}
and
\begin{eqnarray*}
\varepsilon^{-1}\int_0^t\int|Z^IJ|^2wdxd\tau\lesssim\sum\limits_{|K|\leqslant|I|}\varepsilon\int_0^t\int|\partial Z^Kh|^2wdxd\tau
\end{eqnarray*}
\end{lem}
\textbf{Proof.}
The above estimates are straightforward application of Corollary 13.3 of \cite{LR2010} and our Lemma \ref{lem8.2}.
\endproof

Now we deal with $D^I$ and $W^I$.

\begin{lem}\label{lem8.4}
Under the assumptions of Theorem \ref{thm8.1}, we have
\begin{eqnarray}\label{88}
\varepsilon^{-1}\int_0^t\int|D^I|^2wdxd\tau&\lesssim&\varepsilon\sum\limits_{|K|\leqslant|I|}\int_0^t\int\left(|\partial Z^Kh|^2w+|\bar{\partial}Z^Kh|^2w'\right)dxd\tau+\varepsilon^3
\end{eqnarray}
and
\begin{eqnarray}\label{89}
\varepsilon^{-1}\int_0^t\int|W^I|^2wdxd\tau&\lesssim&\varepsilon\sum\limits_{|K|\leqslant|I|}\int_0^t\int\left(|\partial Z^Kh|^2w+|\bar{\partial}Z^Kh|^2w'\right)dxd\tau\\
&&+\varepsilon^3+\varepsilon\sum\limits_{|K|\leqslant|I|}\int_0^t\int|\partial Z^KA|^2w\,dxd\tau\nonumber
\end{eqnarray}
\end{lem}
\textbf{Proof.} We only prove (\ref{89}) since (\ref{88}) follows from the same approach. According to Proposition 5.3 of \cite{LR2010} we arrive at
\begin{eqnarray}\label{90}
&&|\overset{\sim}{\Box}_gZ^IA-\hat{Z}^I\overset{\sim}{\Box}_gA|\\
&\lesssim&\sum\limits_{|K|\leqslant|I|}\sum\limits_{|J|+(|K|-1)_+\leqslant|I|}\left(\frac{|Z^JH|}{1+|q|}+\frac{|Z^JH|_{\mathcal{L}\mathcal{L}}}{1+|q|}\right)|\partial Z^KA|\nonumber\\
&&+\sum\limits_{|K|\leqslant|I|}\left(\sum\limits_{|J|+(|K|-1)_+\leqslant|I|-1}\frac{|Z^JH|_{\mathcal{L}\mathcal{T}}}{1+|q|}+\sum\limits_{|J|+(|K|-1)_+\leqslant|I|-2}\frac{|Z^JH|}{1+|q|}\right)|\partial Z^KA|.\nonumber
\end{eqnarray}
Our goal is to obtain the estimate for the quantity
$$
\sum\limits_{|I|\leqslant N}\int_0^t\int|\overset{\sim}{\Box}_gZ^IA-\hat{Z}^I\overset{\sim}{\Box}_gA|^2w\,dxd\tau.
$$
Let us first deal with the terms in (\ref{90}) with $|K|\leqslant N-4$. In this case we use the decay estimate (\ref{78}). It is clear that now we only have to consider the expression
\begin{equation}\label{91}
\begin{aligned}
&\sum\limits_{\substack{|J|\leqslant|I|\\|J'|\leqslant|I|-1\\|K|\leqslant|I|-2}}\int_0^t\int\left\{\frac{|Z^JH|^2}{(1+|q|)^2}+\frac{|Z^JH|^2_{\mathcal{L}\mathcal{L}}+|Z^{J'}H|^2_{\mathcal{L}\mathcal{T}}+|Z^KH|^2}{(1+|q|)^2}\right\}\varepsilon^2(1+|q|)^{-3}w\,dxd\tau\\
=&\sum\limits_{\substack{|J|\leqslant|I|\\|J'|\leqslant|I|-1\\|K|\leqslant|I|-2}}\int_0^t\int\{|Z^JH|^2+|Z^JH|^2_{\mathcal{L}\mathcal{L}}+|Z^{J'}H|^2_{\mathcal{L}\mathcal{T}}+|Z^KH|^2\}\varepsilon^2(1+|q|)^{-1}w\,dxd\tau.
\end{aligned}
\end{equation}
From the proof of Lemma 11.5 in \cite{LR2010} it follows that (\ref{91}) is bounded by
\begin{eqnarray*}
&&C\varepsilon^2\sum\limits_{|K|\leqslant|I|}\int_0^t\int\left(\frac{|\partial Z^Kh|^2}{1+\tau}w+|\bar{\partial}Z^Kh|^2w'\right)dxd\tau\\
&&+C\varepsilon^2\sum\limits_{|K|\leqslant|I|-1}\int_0^t\int\frac{|\partial Z^Kh|^2}{(1+\tau)^{1-2C\varepsilon}}wdxd\tau+C\varepsilon^4
\end{eqnarray*}
which is equivalent to
\begin{eqnarray*}
&&C\varepsilon^2\sum\limits_{|K|\leqslant|I|}\int_0^t\int\left(|\partial Z^Kh|^2w+|\bar{\partial}Z^Kh|^2w'\right)dxd\tau+C\varepsilon^4,
\end{eqnarray*}
where we let the parameter $M$ in the expression $H_0^{\mu\nu}:=-\chi(r/t)\chi(r)M\delta^{\mu\nu}/r$ equal to 0(the expression can be found at the beginning of the proof of Lemma 11.5 in \cite{LR2010}). For more details we refer to the last inequality on page 1460 of \cite{LR2010}.

Returning to (\ref{90}) we now deal with the case $|K|\geqslant N-3$, which implies $|J|\leqslant4$. By the same way of proving Lemma 11.5 of \cite{LR2010} we arrive at that the contribution of the terms with $|K|\geqslant N-3$ to $|\overset{\sim}{\Box}_gZ^IA-\hat{Z}^I\overset{\sim}{\Box}_gA|$ can be bounded by
$$
\varepsilon\sum\limits_{|K|=|I|}\frac{|\partial Z^KA|}{1+\tau}+\varepsilon\sum\limits_{|K|<|I|}\frac{|\partial Z^KA|}{(1+\tau)^{1-C\varepsilon}},
$$
which is equivalent to
$$
\varepsilon\sum\limits_{|K|\leqslant|I|}|\partial Z^KA|.
$$
\endproof

Now let us finish the proof of Theorem \ref{thm8.1}. Applying (\ref{80}) and (\ref{81}) together with Lemma \ref{lem8.3} and Lemma \ref{lem8.4} yields
\begin{eqnarray}\label{93}
&&\int_{\Sigma_t}(|\partial Z^Ih|^2+|\partial Z^IA|^2)w+\int_0^t\int_{\Sigma_{\tau}}(|\bar{\partial} Z^Ih|^2+|\bar{\partial} Z^IA|^2)w'\\
&\lesssim&\int_{\Sigma_0}(|\partial Z^Ih|^2+|\partial Z^IA|^2)w+\varepsilon\sum\limits_{|K|\leqslant|I|}\int_0^t\int(|\partial Z^Kh|^2+|\partial Z^KA|^2)w\nonumber\\
&&+\varepsilon\sum\limits_{|K|\leqslant|I|}\int_0^t\int(|\bar{\partial} Z^Kh|^2+|\bar{\partial} Z^KA|^2)w'+\varepsilon^3\nonumber.
\end{eqnarray}
Denote
$$\tilde{\mathcal{E}}_k(t):=\sup\limits_{0\leqslant\tau\leqslant t}\sum\limits_{\substack{Z\in\mathfrak{L}\\|I|\leqslant k}}\int_{\Sigma_{\tau}}(|\partial Z^Ih|^2+|\partial Z^IA|^2)w\,dx$$
and
$$
\mathcal{S}_k(t):=\sum\limits_{\substack{Z\in\mathfrak{L}\\|I|\leqslant k}}\int_0^t\int_{\Sigma_{\tau}}(|\bar{\partial} Z^Ih|^2+|\bar{\partial} Z^IA|^2)w'\,dx.
$$
Then we get
\begin{eqnarray}\label{94}
\tilde{\mathcal{E}}_k(t)+\mathcal{S}_k(t)&\lesssim&\tilde{\mathcal{E}}_k(0)+\varepsilon\int_0^t\tilde{\mathcal{E}}_{k}(\tau)d\tau+\varepsilon\mathcal{S}_k(t)+\varepsilon^3,
\end{eqnarray}
which implies
\begin{eqnarray*}
\tilde{\mathcal{E}}_k(t)\leqslant C(T)\varepsilon^2.
\end{eqnarray*}
\endproof

{}

\vspace{1.0cm}

Zonglin Jia

{\small\it Institute of Applied Physics and Computational Mathematics, China Academy of Engineering Physics, Beijing, 100088, P. R. China}

{\small\it Email: 756693084@qq.com}

Boling Guo

{\small\it Institute of Applied Physics and Computational Mathematics, China Academy of Engineering Physics, Beijing, 100088, P. R. China}

{\small\it Email: gbl@iapcm.ac.cn}\\

\end{document}